\documentclass[twoside,11pt,a4paper,leqno]{article}
\usepackage[german,english]{babel}
\usepackage[leqno]{amstex}
\usepackage{amssymb}

\newcommand{\Res}[2]{\mathop{\rm Res}\nolimits_{#1}^{#2}}
\newcommand{\Ind}[2]{\mathop{\rm Ind}\nolimits_{#1}^{#2}}

\newcommand{\Hom}{\mathop{\rm Hom}\nolimits}
\newcommand{\End}{\mathop{\rm End}\nolimits}

\newcommand{\row}{\mathop{\rm row}\nolimits}

\newcommand{\T}{\mbox{{\rm T}}}

\newcommand{\D}{{\cal D}}
\newcommand{\B}{{\cal B}}
\newcommand{\E}{{\cal E}}

\newcommand{\Sch}{{\cal S}}
\newcommand{\sym}[1]{{\frak S}_{#1}}

\newcommand{\Schuh}{$(Q,q)$-Schur algebra}
\newcommand{\pf}{\smallskip\noindent {\bf Proof:\hspace{1em}}}
\let\proof\pf

\def\epf{\ifmmode\eqno\Box\medskip\else{\unskip\nobreak\hfil%
  \penalty50\hskip2em\hbox{}\nobreak\hfil$\Box$
  \parfillskip=0pt \finalhyphendemerits=0\penalty-100\medskip}\fi}
\let\endproof\epf

\renewcommand{\H}{{\cal H}}

\makeatletter
\@addtoreset{equation}{section}
\def\theequation{\thesection.\arabic{equation}}
\makeatother

\newtheorem{thm}[equation]{Theorem}
\newtheorem{lemma}[equation]{Lemma}
\newtheorem{notation}[equation]{Notation}

\newtheorem{cor}[equation]{Corollary}

\newtheorem{rem}[equation]{Remark}
\newtheorem{ex}[equation]{Example}

\newtheorem{ddef}[equation]{Definition}

\newcommand{\mue}[1]{{\mu^{({#1})}}}
\newcommand{\la}[1]{{\lambda^{({#1})}}}
\newcommand{\ta}[1]{\t^{({#1})}}

\newcommand{\tlam}{\t^{\lambda}}

\newcommand{\lahat}{\hat{\lambda}}

\newcommand{\rplus}{\mbox{\boldmath${r}^+$}}
\newcommand{\rminus}{\mbox{\boldmath${r}^-$}}
\renewcommand{\r}{\mbox{\boldmath${r}^\pm$}}
\newcommand{\delo}{{\Delta_0}}

\newcommand{\strich}{\,|\,}

\newcommand{\uap}{u_a^+}
\newcommand{\uam}{u_a^-}
\newcommand{\ubp}{u_b^+}

\newcommand{\bicomp}[1]{\Lambda_2(#1)}
\newcommand{\Mla}{M^\lambda}
\newcommand{\Mmu}{M^\mu}
\newcommand{\toid}[1]{\{#1\}}
\newcommand{\q}{\underline{q}}
\newcommand{\qmlt}[1]{\q^{-\underline{\ell}{(#1)}}}
\newcommand{\qlt}[1]{\q^{\underline{\ell}{(#1)}}}
\newcommand{\Hlm}{{\frak H}_{\lambda,\mu}}
\newcommand{\dou}{\D_{{\tilde{B}},{\tilde{A}}}}
\newcommand{\tA}{{\tilde{A}}}
\newcommand{\tB}{{\tilde{B}}}
\newcommand{\tBd}{{\tilde{B}d}}
\newcommand{\tBdA}{{\tilde{B}d\cap\tilde{A}}}
\newcommand{\BdA}{{Bd\cap A}}
\newcommand{\dHlm}{{\frak H}_{\lambda,\mu}^{(d)}}
\newcommand{\vdHlm}{{\frak H}_{\lambda,\mu}^{(d,v)}}
\newcommand{\uvdHlm}{{\frak H}_{\lambda,\mu}^{(d,v,u)}}
\newcommand{\lmphi}{\varphi_{\mu,\lambda}^{(d,v,u)}}
\newcommand{\lmcphi}{\varphi_{\mu,\lambda}^c}
\newcommand{\IJhatphi}{\hat{\varphi}_{I,J}^{(d,v,u)}}
\newcommand{\lmhatphi}{\hat{\varphi}_{\mu,\lambda}^{(d,v,u)}}
\newcommand{\llphi}{\varphi_{\lambda,\lambda}}
\newcommand{\lophi}{\varphi_{\lambda,\om}}
\newcommand{\oophi}{\varphi_{\om,\om}}
\newcommand{\stB}{{\frak B}}
\newcommand{\stBlm}{{\frak B}_{\lambda,\mu}}
\newcommand{\zubar}{\underline{z}}
\newcommand{\IJtphi}{\tilde{\varphi}_{I,J}^c}
\newcommand{\om}{\omega}
\newcommand{\btabl}{{(\T^{(1)},\T^{(2)})}}
\newcommand{\Ta}[1]{\T^{(#1)}}
\newcommand{\Bit}{\mathfrak{T}(\lambda,\mu)}
\newcommand{\Bittilde}{\mathfrak{T}(\tilde{\lambda},\tilde{\mu})}
\newcommand{\Bitplus}{\mathfrak{T}^+(\lambda,\mu)}
\newcommand{\Bitss}{\mathfrak{T}_\circ(\lambda,\mu)}
\newcommand{\ep}{\varepsilon}
\newcommand{\tilmu}[1]{\tilde{\mu}_{#1}}

\let\len\ell
\makeatletter

\def\enumerate{%
  \ifnum \@enumdepth >3 \@toodeep
  \else\advance\@enumdepth \@ne
       \edef\@enumctr{enum\romannumeral\the\@enumdepth}%
       \topsep\z@\parskip\z@\partopsep\z@                                     
       \list
       {\csname label\@enumctr\endcsname}
       {\@nmbrlisttrue\setcounter{\@enumctr}{0}\let\@listctr\@enumctr
       \def\makelabel##1{\hss\llap{##1}}
       \parsep\z@\itemsep\z@}%
  \fi}

\def\Number#1{\refstepcounter{equation}%
              \leqno(\theequation)\if*#1\else\label{#1}\fi}     

\newtheorem{res}[equation]{\@gobble}

\makeatother

\let\Medskip\relax
\def\s{\mathfrak s}
\def\t{\mathfrak t}
\def\u{\mathfrak u}

\def\T{\mathtt T}
\def\A{\mathtt A}
\def\B{\mathtt B}

\def\({\Big(}
\def\){\Big)}
\let\<\langle
\let\>\rangle
\def\Ll{{\ll}}
\def\Gg{{\gg}}

\def\Sum{\displaystyle\sum}

\begin{document}
\thispagestyle{empty}

\title{The $(Q,q)$--Schur Algebra}
\author{
            Richard Dipper\\
            Mathematische Institut B\\
            Universit{\"a}t Stuttgart\\
            Postfach 80 11 40\\
            70550 Stuttgart \\
            Deutschland
\and        
            Gordon James\\
            Department of Mathematics\\
            Imperial College\\
            Queen's Gate\\
            London SW7 2BZ\\
            England
\and
            Andrew Mathas\\
            Department of Mathematics\\
            Imperial College\\
            Queen's Gate\\
            London SW7 2BZ\\
            England
}
\makeatletter
\let\@makefnmark\relax
\footnotetext{{\it \kern-6mm%
A.M.S. subject classification (1991): 16G99, 20C20, 20G05}\\
This paper is a contribution to the DFG project on ``Algorithmic
number theory and algebra''. The authors acknowledge support
from DFG; the third author was also supported in part by
SERC grant GR/J37690}
\makeatother

\maketitle

\begin{abstract}
In this paper we use the Hecke algebra of type $\bf B$ to define a new 
algebra $\Sch$ which is an analogue of the $q$--Schur algebra. We construct 
Weyl modules for $\Sch$ and obtain, as factor modules, a family of
irreducible $\Sch$--modules over any field.
\end{abstract}

\section{Introduction}

The ordinary Schur algebra is of key importance in the study of the
representation theory of general linear groups in the describing
characteristic, and it provides a link between the general linear groups and 
the symmetric groups. In \cite{DJ4,DJ5} we introduced the $q$--Schur
algebra, and demonstrated its usefulness in the representation theory
of $GL_n(q)$ over a field of non--describing characteristic. In 
\cite{DDo} it was shown that the $q$--Schur algebra is given as the dual
of a homogeneous part of quantum--$GL_n$, or alternatively as the
factor of quantum--$GL_n$ or the corresponding quantum enveloping
algebra modulo the kernel of its action on quantum tensor space
(compare \cite{BLM}). In particular, the
representations of $q$--Schur algebras are precisely the homogeneous 
polynomial representations of quantum--$GL_n$ in a fixed degree
(compare \cite {D7}). 

The construction of the $q$--Schur algebra involves the Hecke algebra of type 
$\bf A$; in this paper we use the Hecke algebra of type $\bf B$ to
build an algebra which we call the $(Q,q)$--Schur algebra. Others have
devised a version of a Schur algebra of type $\bf B$ \cite{GH,G}, but ours
is a larger algebra. Applications to the representation
theory of finite symplectic groups in the non--describing
characteristic case have already been provided \cite{Gru,DG}, and we
expect that further applications will ensue, using our larger and more
complicated algebra.

Hecke algebras of type $\bf B$ have been studied in \cite{DJ6,DJM1,M:can}. We
begin by recalling and extending some of the notation and results which we 
used in those papers. The remainder of the paper is then devoted to
introducing the $(Q,q)$--Schur algebra  and investigating its main
properties. 
In particular we construct a generic basis of the $(Q,q)$--Schur
algebra and we define
$(Q,q)$--Weyl modules. The Weyl modules are
labelled by bipartitions and they have unique maximal submodules. The
corresponding factor modules are pairwise non--isomorphic
irreducible representations of the $(Q,q)$--Schur algebra. We show that the
decomposition matrix
which describes the composition multiplicities of these irreducible
modules in the $(Q,q)$--Weyl modules is unitriangular, and we construct a 
``semistandard basis'' for each $(Q,q)$--Weyl module.

In a forthcoming paper we shall construct a cellular basis of the \Schuh\ 
$\Sch$.
As a consequence, every irreducible representation of $\Sch$ is isomorphic
to one of the irreducible representations which we construct here and, in
addition, $\Sch$ is quasi--hereditary. Furthermore, we shall generalize some
of our 
results to construct Schur algebras of the Ariki--Koike algebras.


\section{The Hecke algebra of type $B$}

Let $W_r$ be the group $C_2\wr\sym{r}$, where $\sym{r}$ is the
symmetric group of degree $r$. Then $W_r$ is generated by elements
$s_0,s_1,\ldots , s_{r-1}$ which satisfy the following relations:
$$\begin{array}{rcl@{\quad}l}
s_i^2 & = & 1 & \text{for $0\leq i\leq r-1$}\\
s_is_j & = & s_js_i & \text{if $1\leq i+1< j\leq r-1$}\\
s_is_{i+1}s_i & =  &s_{i+1}s_is_{i+1} & \text{if $1\leq i\leq r-2$}\\
s_0s_1s_0s_1 & =  &s_1s_0s_1s_0.
\end{array}\Number{2-1}$$

Let $\rplus = \{1,2,\ldots ,r\}$ and $\rminus = \{-1,-2,\ldots
,-r\}$. We identify $W_r$ with a subgroup of the symmetric group on
$\r = \rplus \cup \rminus$ by letting
\begin{align*}
s_0 & = (1,-1)\\
s_i & = (i,i+1)(-i,-i-1) \qquad \text{for $1\leq i\leq r-1$}.
\end{align*}
It is easy to see that $W_r$ acts transitively on $\r$. The subgroup 
of $W_r$ generated by $s_1,\ldots ,s_{r-1}$ is the symmetric group 
$\sym{r}$ of degree $r$.

The concepts of a {\bf reduced expression} for an element of $W_r$ and
the {\bf length} $\len(w)$ of $w\in W_r$ are defined in the usual way. It
is useful to note the following well-known method for calculating
$\len(w)$. Consider the ${\Bbb Q}W_r$--module $V$ whose basis is $\{e_i\,|\,
i\in \rplus\}$ and where, for $i\in\rplus$ and $w\in W_r$, we have 
\begin{equation*}
e_iw = \begin{cases}
e_{iw} & \text{if $iw > 0$} \\
-e_{-iw}& \text{if $iw < 0$}.
\end{cases}
\end{equation*}
 The {\bf roots}
for $W_r$ are 
$$
\pm e_i \  \text{and} \  \pm e_i\pm e_j \quad (i,j \in
\rplus, i \neq j).
$$
The {\bf positive roots} are 
$$
e_i,\  e_i+e_j,\  e_j-e_i \quad (i,j\in \rplus, i < j),
$$
and all other roots are {\bf negative roots}.
The length of $w$ is equal to the number of positive roots changed to
negative roots by $w$.

The simple roots
$$\alpha_0=e_1,\;\alpha_1=e_2-e_1,\;\alpha_2=e_3-e_2,\;\ldots,\;
         \alpha_{r-1}=e_r-e_{r-1}$$
form a basis $\Delta_0$ of $V$, and every positive root is a non--negative
linear combination of simple roots. Moreover, for $i=0,1,\ldots,r-1$, $s_i$
acts on $V$ as the reflection in the hyperplane orthogonal to $\alpha_i$.
If $\alpha$ is a positive root we sometimes write $s_\alpha$ for the 
reflection in the hyperplane orthogonal to $\alpha$; in particular,
$s_i=s_{\alpha_i}$. Note that $s_\alpha\in W_r$ for all positive roots
$\alpha$.

Let $R$ be a commutative ring with $1$, and let $q$ and $Q$ be invertible
elements of $R$. The {\bf Hecke algebra} $\H_{R,q,Q}(W_r)$ {\bf of
  type} $B_r$ is defined to be the free $R$--module with basis
$\{T_w\strich w\in W_r\}$ and multiplication defined as below. If $e$ is
the identity element of $W$ then $T_e$ is the multiplicative identity for 
$\H_{R,q,Q}(W_r)$ and for $\rho\in R$ we abbreviate $\rho T_e$ as $\rho$. We 
often write $T_{s_i}$ as $T_i$ for $0\leq i\leq r-1$ and $\H_{R,q,Q}(W_r)$ as 
$\H (W_r)$ or as $\H$. Then 
\begin{enumerate}
\item if $w = v_1v_2\cdots v_l$ is a reduced expression for $w\in W_r$
  where each $v_i$ belongs to $\{s_0,s_1,\ldots ,s_{r-1}\}$, then $T_w
  = T_{v_1}T_{v_2}\cdots T_{v_l}$;
\item $(T_i)^2 = q + (q-1)T_i$ for $1\leq i\leq r-1$;
\item $(T_0)^2 = Q + (Q-1)T_0$.
\end{enumerate}

Let $\H(\sym{r})$ denote the subalgebra of $\H$ generated by
$T_1,T_2,\ldots ,T_{r-1}$. 

For each pair of integers $i,j$ in $\rplus$ define $s_{i,j}\in W_r$
by
$$
s_{i,j} = \left\{\aligned s_is_{i+1}\ldots s_{j-1} &\qquad \text{if $i\le j$} \\
                          s_{i-1}s_{i-2}\ldots s_j &\qquad \text{if $i> j$}
                 \endaligned \right.
$$  

For example, $s_{1,r}$ is the permutation of $\r$ which is given
in its cycle decomposition as
$$
(r,r-1,\ldots ,2,1)(-r,-r+1,\ldots ,-2,-1).
$$
For $0\leq a\leq r$, let $w_{a,r-a} = (s_{1,r})^a$. As in
\cite[2.4--2.9]{DJ6} we have the following.

\begin{res} Let $i,j\in \rplus$ and $0\leq a\leq r$. Then
\begin{enumerate}
\item $\len(s_{i,j}) = |j-i|$.
\item $\len(w_{a,r-a}) = a(r-a)$.
\item $s_{a+1,1}s_{a+2,2}\cdots s_{r,r-a}$ gives a reduced expression
  for $w_{a,r-a}$.
\end{enumerate}
\label{2.1}\end{res}

For typographical reasons, let $T_{i,j} = T_{s_{i,j}}$ and $h_{a,r-a}
= T_{w_{a,r-a}}$. 

\begin{ddef} For $0\leq a\leq r$, let the elements $\uap$ and $\uam$
of $\H$ be given~by
\begin{align*}
\uap &= \prod_{i=1}^a(q^{i-1}+T_{i,1}T_0T_{1,i})\\
\uam &= \prod_{i=1}^a(Qq^{i-1}-T_{i,1}T_0T_{1,i})
\end{align*}
\end{ddef}

It is easy to prove the following (compare \cite[3.3 and 3.4]{DJ6}).
\begin{res} Let $0\leq a\leq r$.  
\begin{enumerate}
\item If $0\leq b\leq r-1$ with $b\neq a$, then $\uap$ and $\uam$ 
      commute with $T_b$.
\item If $a>0$ then $\uap T_0 = Q\uap$ and $\uam T_0 = -\uam$.
\end{enumerate}\label{2.3}
\end{res}

An $a$-{\bf bicomposition} of $r$ is an ordered pair $(\la{1},\la{2})$
of compositions, where $\la{1}$ is a composition of $a$ and $\la{2}$
is a composition of $r\!-\!a$; if both $\la{1}$ and $\la{2}$ are
partitions, then $(\la{1},\la{2})$ is an $a$-{\bf bipartition} of
$r$. 

Let $\bicomp{n,r}$ denote the set of bicompositions $\lambda =
(\la{1},\la{2})$ of $r$ with the property that the sum of the number
of parts of $\la{1}$ and the number of parts of $\la{2}$ is $n$. (Note
that we allow a composition to have zero parts.)

We order the bipartitions of $r$ by letting all $a$-bipartitions
precede all $b$-bipartitions if $a>b$, and by saying that the $a$-bipartition
$\lambda = (\la{1},\la{2})$ precedes the $a$-bipartition $\mu =
(\mue{1},\mue{2})$ if $\la{1}$ precedes $\mue{1}$ lexicographically or
$\la{1} = \mue{1}$ and $\la{2}$ precedes $\mue{2}$
lexicographically. We call this the {\bf lexicographic order} on
bipartitions.

Suppose that $\lambda = (\la{1},\la{2})$ and $\mu = (\mue{1},\mue{2})$
are bicompositions of~$r$. We say that $\lambda$ and $\mu$ are {\bf 
associated} if $\mue{i}$ can be obtained from $\la{i}$, $i=1,2$, by
reordering the parts.

>From $\lambda$ we obtain a corresponding {\bf diagram} $[\lambda]$ which
consists of crosses in the plane, as in the example $\lambda =
((4,3,1),(3,2))$, where 
 
$$
[\lambda] = \left(\begin{matrix} \times & \times & \times & \times \\
                                 \times & \times & \times & \      \\
                                 \times & \      & \      & \   \end{matrix}
                     \quad , \quad
                  \begin{matrix} \times & \times & \times\\
                                 \times & \times & \   \\
                                 \      & \      & \    \end{matrix}
             \right )
$$

A $\lambda$-{\bf bitableau} is obtained from $[\lambda]$ by replacing
each cross by one of the numbers from $\r$ in such a way that for
every $i\in \rplus,$ precisely one of $\{i,-i\}$ is used. If $\t$ is a
$\lambda$--bitableau then $\t = (\ta{1},\ta{2})$, where $\ta{1}$ is a
$\la{1}$-tableau and $\ta{2}$ is a $\la{2}$-tableau.

\begin{ddef} Suppose that $\lambda= (\la{1},\la{2})\in\bicomp{n,r}$,
  and that $\t = (\ta{1},\ta{2})$ is a $\lambda$--bitableau. Assume that
  $i\in\t$. Let $\row_\t(i) = j$ if $i$ belongs to row $j$ of $\ta{1}$
  and let $\row_\t(i) = n+j$ if $i$ belongs to row $j$ of $\ta{2}$.
\end{ddef}

\begin{ddef} Let $\lambda$ be a bicomposition.
\begin{enumerate}
\item A $\lambda$--bitableau $\t = (\ta{1},\ta{2})$ is {\bf row 
standard} if the entries increase from left to right in each row 
of $\ta{1}$ and in each row of $\ta{2}$, and all entries in 
$\ta{1}$ belong to $\rplus$.
\item The $\lambda$--bitableaux $\t = (\ta{1},\ta{2})$ and 
  $\s= (\s^{(1)},\s^{(2)})$ are 
  {\bf row equivalent} if $|\ta{1}|$ and
  $|\s^{(1)}|$ are row equivalent, and $\ta{2}$ and $\s^{(2)}$ are
  row-equivalent.
  Here the entries of $|\ta{1}|$ are the absolute values of the
  entries of $\ta{1}$ and $|\s^{(1)}|$ is defined similarly. A row
  equivalence class of the $\lambda$--bitableaux is a
  {\bf $\lambda$--bitabloid} and the $\lambda$-bitabloid containing $\t$ is
  denoted by $\{\t\}$. 
\item A $\lambda$--bitableau is {\bf standard} if all its entries
belong to $\rplus$, and it is row standard, and all the entries
increase down each column of $\ta{1}$ and each column of $\ta{2}$.
\end{enumerate}
\label{2.5}\end{ddef}

We remark that every $\lambda$--bitabloid contains exactly one row standard
$\lambda$--bitableau. When dealing with $\lambda$-bitabloids
$\toid{\t}$, we often find it convenient to specify that $\t =
(\ta{1},\ta{2})$ is row standard; this ensures that all the entries in
$\ta{1}$ are positive.

Here is an example of a row standard bitableau:

$$
 \left(\begin{matrix} 2 & 8 \\
                      3 & 4 
                      \end{matrix}
                     \quad, \quad 
                  \begin{matrix} -5 & -1 & 6 \\
                                 -7 &\ 9 & \  \end{matrix}
             \right).
$$  

Next, we wish to specify, for a given $\lambda = (\la{1},\la{2})$, 
several special standard $\lambda$--bitableaux. The 
$\lambda$--bitableaux which we wish to define are most easily
understood with an example.

\begin{ex} Suppose that $\lambda = ((4,3,1),(3,2))$. We shall define
  the $\lambda$--bitableaux $\t^\lambda$, $\hat{\t}^\lambda$, $\t_\lambda$
  and $\hat{\t}_\lambda$ so that

\begin{alignat*}{4}
\t^\lambda & = & &\left(\begin{matrix}         \ 1 &\ 2 &\ 3 &\ 4 \\
                                             \ 5 &\ 6 &\ 7 &\ \ \\
                                             \ 8 &\ \ &\ \ &\ \   
\end{matrix}\right.
                   \quad &,&& \quad 
                  \left.\begin{matrix}       \ 9 & 10 & 11 \\
                                              12 & 13 & \ \ \\
                                             \ \ & \  & \ \  \end{matrix}
\right)&,\\
\hat{\t}^\lambda  & = & &\left(\begin{matrix}   \ 6 &\ 7 &\ 8 &\ 9 \\
                                              10 & 11 & 12 &\ \ \\
                                              13 &\ \ &\ \ & \ \   
\end{matrix}\right.
                    \quad &,&& \quad 
                  \left.\begin{matrix}       \ 1 &\ 2 &\ 3 \\
                                             \ 4 &\ 5 &\ \ \\
                                             \ \ &\ \ & \ \ \end{matrix}
            \right)&,\\
\t_\lambda  & =& &\left(\begin{matrix}         \ 1 &\ 4 &\ 6 &\ 8 \\
                                             \ 2 &\ 5 &\ 7 &\ \ \\
                                             \ 3 &\ \ &\ \ &\ \  
\end{matrix}\right.
                     \quad &,&&\quad 
                  \left.\begin{matrix}      \ 9  & 11 & 13 \\
                                              10 & 12 &\ \ \\
                                             \ \ &\ \ &\ \  \end{matrix}
             \right)&,\\
\hat{\t}_\lambda & =& &\left(\begin{matrix}    \ 6 &\ 9 & 11 & 13 \\
                                             \ 7 & 10 & 12 &\ \ \\
                                             \ 8 &\ \ &\ \ &\ \  
\end{matrix}\right.
                     \quad &,&& \quad 
                  \left.\begin{matrix}       \ 1 &\ 3 &\ 5 \\
                                             \ 2 &\ 4 &\ \ \\
                                             \ \ &\ \ &\ \ \end{matrix}
             \right)&.
\end{alignat*}
\end{ex}
\Medskip

\begin{ddef}\label{2.7}
Suppose that $\lambda$ is an $a$-bicomposition of $r$.
\begin{enumerate}
\item Let $\t^\lambda$ = $(\t^{\la{1}},\t^{\la{2}})$ be the standard
  $\lambda$--bitableau in which the numbers $1,2,\ldots, a$ appear in
  order by rows in $\t^{\la{1}}$ and the numbers $a+1,a+2,\ldots ,r$
  appear in order by rows in $\t^{\la{2}}$.  
\item Let $\hat{\t}^\lambda = (\hat{\t}^{\la{1}},\hat{\t}^{\la{2}})$ be
  the standard $\lambda$--bitableau in which the numbers $1,2,\ldots
  ,r-a$ appear in order by rows in $\hat{\t}^{\la{2}}$ and the numbers
  $r-a+1,r-a+2,\ldots ,r$ appear in order by rows in
  $\hat{\t}^{\la{1}}$.
\item Let $\t_\lambda = ({\t_\lambda}^{(1)},{\t_\lambda}^{(2)})$ be the standard
  $\lambda$--bitableau in which the numbers $1,2,\ldots ,a$ appear in
  order by columns in ${\t_\lambda}^{(1)}$ and the numbers $a+1,a+2,\ldots ,r$
  appear in order by columns in ${\t_\lambda}^{(2)}$. 
\item Let 
  $\hat{\t}_\lambda=({{\hat{\t}}_\lambda}^{\phantom{\lambda}(1)},
                     {{\hat{\t}}_\lambda}^{\phantom{\lambda}(2)})$ 
  be the standard $\lambda$--bitableau in which the numbers $1,2,\ldots ,r-a$
  appear in order by columns in 
  ${{\hat{\t}_\lambda}}^{\phantom{\lambda}(2)}$ and 
  the numbers $r-a+1,r-a+2,\ldots ,r$ appear in order by columns in 
  ${{\hat{\t}}_\lambda}^{\phantom{\lambda}(1)}$.
\end{enumerate}
\end{ddef}
\Medskip

Note that $W_r$ acts on the set of $\lambda$--bitableaux. If $\t$ is a
$\lambda$--bitableau and $w\in W_r$, then we obtain the
$\lambda$--bitableau $\t w$ by replacing each $i$ in $\t$ by $iw$. For
example, if $\lambda$ is an $a$-bicomposition of $r$, then $\t^\lambda
w_{a,r-a} = \hat{\t}^\lambda$. It is easy to see that $W_r$ acts 
transitively on the set of $\lambda$--bitableaux. 

One easily checks that the action of $W_r$ preserves the row equivalence
classes of $\lambda$--bitableaux. We therefore have a transitive action of 
$W_r$ on the set of $\lambda$--bitabloids. If
$\lambda$ is an $a$-bicomposition, then the stabilizer of the
bitabloid $\{\t^\lambda\}$ is the subgroup
$$
W_\lambda = ((\underbrace{C_2\times\cdots\times C_2}_{\text{$a$
factors}})\rtimes{\frak S}_\la{1}) \times{\frak S}_\la{2}
\Number{W-lambda}$$
of $W_r$. 
Thus the permutation representation of $W_r$ on the set of
$\lambda$--bitabloids is equivalent to the representation of $W_r$ on the
cosets of $W_\lambda$.

A {\bf reflection subgroup} $W_J$ of $W_r$ is a subgroup generated by a set of
reflections $\{s_\alpha\,|\,\alpha\in J\}$ for some subset $J$ of the positive
roots. In particular, if $J\subseteq\delo$ then $W_J$ is a {\bf parabolic
subgroup} of $W_r$. If $\lambda$ is an $a$--bicomposition then
$W_\lambda$ is a parabolic subgroup of $W_r$ if and only if 
$\la{1}=(a)$; in this case $W_\lambda$ is a direct product of a Weyl group 
$W_a$ of type~$\bf B$ and several Weyl groups of type $\bf A$. In
general the reflection group $W_\lambda$ has precisely one factor of 
type $\bf B$ for each nonzero part of $\la{1}$ .

Associated with the parabolic subgroup $W_J$, $J\subseteq\delo$, is the
parabolic subalgebra $\H_J=\sum_{w\in W_J}R T_w$ of $\H$. This algebra is 
isomorphic to the Hecke algebra $\H_{R,Q,q}(W_J)$ of $W_J$. If $W_\lambda$ 
is a reflection subgroup of $W_r$, which is not a parabolic subgroup, then
$\sum_{w\in W_\lambda}R T_w$ is not a subalgebra of $\H$ in general.

\begin{ddef} Suppose that $\lambda\in\bicomp{n,r}$. The elements
  $x_\lambda$, $\hat{x}_\lambda$, $y_\lambda$ and $\hat{y}_\lambda$ of
  $\H$ are defined as follows.
\begin{align*}
x_\lambda &= \sum\{T_w\strich w\in W_r \text{\ and $w$ stabilizes the
  rows of $\t^\lambda$}\}\\
\hat{x}_\lambda &= \sum\{T_w\strich w\in W_r \text{\ and $w$ stabilizes the
  rows of $\hat{\t}^\lambda$}\}\\
y_\lambda &= \sum\{(-q)^{-\len(w)}T_w
\strich w\in W_r \text{\ and $w$ stabilizes the
  columns of $\t_\lambda$}\}\\
\hat{y}_\lambda &= \sum\{(-q)^{-\len(w)}T_w
\strich w\in W_r \text{\ and $w$ stabilizes the
  columns of $\hat{\t}_\lambda$}\}
\end{align*}
\label{x-lambda defn}\end{ddef}
\Medskip

\begin{ex} Suppose that $\lambda = ((2,1),(1))$. Then
\begin{xalignat*}{2}
\t^\lambda &=  \left(\begin{matrix} 1 & 2 \\
                                   3 & \  
                      \end{matrix}
                     \quad , \quad 
                  \begin{matrix}  4   \\
                                  \   \end{matrix}
             \right) & \qquad
\t_\lambda &= \left(\begin{matrix} 1 & 3 \\
                                  2 & \ \end{matrix}
                                  \quad , \quad 
                    \begin{matrix} 4 \\
                                   \   \end{matrix} \right) \notag \\
\hat{\t}^\lambda &= \left(\begin{matrix} 2 & 3 \\
                                        4 & \ 
                         \end{matrix}
                         \quad , \quad  
                         \begin{matrix} 1  \\
                                       \    \end{matrix}\right) & \qquad
\hat{\t}_\lambda &= \left(\begin{matrix} 2 & 4 \\
                                  3 & \ \end{matrix}
                                  \quad , \quad 
                    \begin{matrix} 1 \\
                                   \  \end{matrix} \right)\notag 
\end{xalignat*}
and 
\begin{xalignat*}{2}
x_\lambda &= 1 + T_{(1,2)} & \qquad y_\lambda &= 1 -
q^{-1}T_{(1,2)}\notag \\
\hat{x}_\lambda &= 1+T_{(2,3)} &\qquad \hat{y}_\lambda &=
1-q^{-1}T_{(2,3)}.\notag
\end{xalignat*}
\end{ex}

We remark that all four elements defined in \eqref{x-lambda defn}
belong to $\H(\sym{r})$ and  that
$ \hat{x}_\lambda = x_{\hat{\lambda}}$ and $\hat{y}_\lambda=y_{\hat{\lambda}}$,
where $\lahat$ is the $(r-a)$-bicomposition $(\la{2},\la{1})$. 

>From (\ref{2.3}) we obtain the following.
\begin{res} Assume that $\lambda$ is an $a$-bicomposition. Then
\begin{enumerate}
\item $x_\lambda$ and $y_\lambda$ commute with $\uap$;
\item $\hat{x}_\lambda$ and $\hat{y}_\lambda$ commute with
  $u_{r-a}^-$.
\end{enumerate}
\label{2.10}\end{res}

As we noted in \cite[(2.5)]{DJ6}, we also have the following (see the
remark after (\ref{2.1})).

\begin{res} If $\lambda$ is an $a$-bicomposition then
$
x_\lambda h_{a,r-a} = h_{a,r-a}\hat{x}_\lambda.
$
\label{2.11}\end{res}

\begin{ddef} Let $\pi_\lambda$ be the element of $W_r$ which is given
  by $\t^\lambda\pi_\lambda =\t_\lambda$ and $\hat{\pi}_\lambda$ be the
  element of $W_r$ which is given by 
$\hat{\t}^\lambda\hat{\pi}_\lambda = \hat{\t}_\lambda$
\label{2.12}\end{ddef}

We note that $\pi_\lambda$ and $\hat{\pi}_\lambda$ belong to
$\sym{r}\leq W_r$, since the entries of the bitableaux involved are all
positive. Moreover $\hat{\pi}_\lambda = \pi_{\hat{\lambda}}$, and
$h_{a,r-a}T_{\hat{\pi}_\lambda}=T_{\pi_\lambda}h_{a,r-a}$. Note too that
$$
T_{\hat{\pi}_\lambda}u_{r-a}^- = u_{r-a}^-T_{\hat{\pi}_\lambda}
$$ 
by \eqref{2.3} since $\hat{\pi}_\lambda\in\sym{(r-a,a)}$.

Now, \cite[(3.11)]{DJ6} and \cite[(4.1)]{DJ1} imply the following
fundamental result.

\begin{thm} Suppose that $\lambda$ is an $a$-bipartition. Then
$$
\uap x_\lambda \H u_{r-a}^- \hat{y}_\lambda
$$
is a one-dimensional $R$--module. It is spanned by the vector
$$
\uap x_\lambda h_{a,r-a}
T_{\hat{\pi}_\lambda}u_{r-a}^-\hat{y}_\lambda.
$$
\label{2.13}\end{thm}
\medskip

\begin{ddef} Suppose that $\lambda$ is an $a$-bipartition. Let
$$
z_\lambda = \uap x_\lambda h_{a,r-a}T_{\hat{\pi}_\lambda}
  u_{r-a}^-\hat{y}_\lambda.
$$
\label{z-lambda def}\end{ddef}

Note that, in the light of (\ref{2.10}) and (\ref{2.11}), we have several
alternative expressions for $z_\lambda$; for example,
$$\begin{array}{rl}
z_\lambda &=\uap x_\lambda T_{\pi_\lambda}h_{a,r-a} u_{r-a}^-\hat{y}_\lambda\\
          &= \uap h_{a,r-a}u_{r-a}^- \hat{x}_\lambda
                T_{\hat{\pi}_\lambda}\hat{y}_\lambda\\
          &= \uap x_\lambda h_{a,r-a}u_{r-a}^-
                T_{\hat{\pi}_\lambda}\hat{y}_\lambda.
\end{array}\Number{2-3}$$

We call the right ideal $z_\lambda\H$ of $\H$ the $(Q,q)$--{\bf Specht
module} $S^\lambda$. It is straightforward to show that $S^\lambda$
is the dual of the module $\tilde{S}^\lambda$ which we introduced in
\cite[\S 4]{DJM1}; in particular $S^\lambda$ has dimension equal to
the number of standard $\lambda$--bitableaux. Moreover 
there exists an $\H$--invariant bilinear form 
on $M^\lambda = \uap x_\lambda\H$  and a submodule theorem holds for 
$M^\lambda$ (compare Theorem~\ref{submodulethm} and \cite[4.8]{DJ1}).
When $R$ is a field, the submodule theorem allows us to construct
irreducible $\H$--modules, and so we can recover the main results of 
\cite{DJM1}. We do not wish to
pursue this here; instead, we turn to the construction of the
$(Q,q)$-Schur algebra.
  


\section{The $\H$--module $M^\lambda$}

\begin{ddef} Suppose that $\lambda$ is an $a$-bicomposition . Let
  $M^\lambda$ be the right ideal of $\H$ which is given by 
$ M^\lambda = \uap x_\lambda\H$.
\label{3.1}\end{ddef}
\Medskip

If $\lambda$ is an $a$--bicomposition then $M^\lambda$ contains the 
$(Q,q)$-Specht module~$S^\lambda$,
and postmultiplication by $u^-_{r-a}\hat{y}_\lambda$ maps $\Mla$ on
to the one-dimensional $R$--module spanned by the generator $z_\lambda$ of
$S^\lambda$, by Theorem \ref{2.13}.

An argument similar to that of \cite[(1.1)]{DJ5} shows the following.

\begin{res} If $\lambda$ and $\mu$ are associated bicompositions of
  $r$, then $\Mla$ and $\Mmu$ are isomorphic $\H$--modules.
\label{3.2}\end{res}
\Medskip

We next construct a basis of $\Mla$.

\begin{thm} Assume that $\lambda$ is an $a$-bicomposition of $r$. Then
$\Mla$ is a free $R$--module with basis
$$
\{\uap x_\lambda T_w\strich w\in W_r \text{\ and $\t^\lambda w$ is row
  standard}\}.
$$
\label{3.3}\end{thm} 
\pf Note first, that 
$ \uap x_\lambda\H(W_a) = \uap x_\lambda\H(\sym{a})$.
This is immediate if $a=0$ since $\H(W_0) = \H(\sym{0}) = \{1\}$. If
$a>0$ then from \eqref{2.3}, for $h\in\H(\sym a)$,
\begin{equation*}
\uap x_\lambda h T_0 = x_\lambda h\uap T_0 = Qx_\lambda h\uap = Q\uap
x_\lambda h.
\end{equation*}
Next, \cite[(3.2)(i)]{DJ1} implies that a basis of
$\uap x_\lambda\H(\sym{a})$ is given by 
\begin{equation*}
\{\uap x_\lambda T_w\strich w\in\sym{a} \text{\ and $\t^\lambda w$ is
  row standard}\}.
\end{equation*}
Let ${\frak S}_\la{2}$ denote the subgroup of ${\frak
  S}_{\{a+1,\ldots,r\}}$ which stabilizes the rows of $\t^\la{2}$. Then
$W_a\times{\frak S}_\la{2}$ is a parabolic subgroup of $W_r$. Let $\D$
denote the set of distinguished right coset representatives of
$W_a\times{\frak S}_\la{2}$ in $W_r$. We claim that $\B$ is a basis of
$\Mla$, where
$$\B = \{\uap x_\lambda T_wT_d\strich w\in{\frak S}_a, d\in\D \text{\
  and $\t^\lambda w$ is row standard}\}.
\Number{3-3}$$
To see this notice that we
can write any element in $W_r$ as $w_1w_2d$ where $w_1\in W_a$,
$w_2\in {\frak S}_{\la{2}}$ and $d\in\D$; but 
$$
\uap x_\lambda T_{w_1}T_{w_2}T_d = q^{\len(w_2)}\uap x_\lambda T_{w_1}T_d
\in \uap x_\lambda\H(\sym{a})T_d,
$$
so $\B$ spans $M^\lambda$. Secondly, if $d$ and $d'$ are distinct
elements of $\D$ then corresponding elements in $\B$ have distinct
supports, and hence $\B$ is linearly independent.

Since $T_wT_d = T_{wd}$ for $w\in\sym{a}$ and $d\in\D$, the proof of
the theorem will be complete if we show that
\begin{multline}
\{wd\strich w\in\sym{a},\  d\in\D \text{\ and $\t^\lambda w$ is
  row standard}\}\\ = \{v\in W_r\strich \t^\lambda v \text{\ is
  row standard}\}.
\label{3-4}\end{multline}
Suppose that $w\in\sym{a}$, $d\in\D$ and $\t^\lambda w$ is row standard.
Since $d$ is a distinguished coset representative, $d$ sends the
positive roots for $W_a\times{\frak S}_\la{2}$ into positive
roots. The positive roots for $W_a\times{\frak S}_\la{2}$ are
$$
e_i\quad (1\leq i\leq a),\ e_i+e_j,\ e_j-e_i\quad (1\leq i<j\leq a), 
$$
together with
$$ 
e_j-e_i \quad  (i<j \text{\ with $i$, $j$ in the same row of\ } \t^\la{2}).
$$
If $1\leq i\leq a$ then $e_i d$ is a positive root, so $id\in\rplus$; thus
all the entries in the first component of $\t^\lambda wd$ must be positive. 
If $i<j$ and $i$ and $j$ are in the same row of $\t^\lambda w$ then 
$(e_j-e_i)d$ is a positive root, so $id<jd$. Therefore, $\t^\lambda wd$ is
row standard (see Definition \ref{2.5}).

Finally, note that the sets on the two sides of equation \eqref{3-4} have
the same size. This completes the proof of \eqref{3-4} and hence of
the theorem.
\epf

\begin{ddef} Let $J\subseteq \delo$. The subalgebra of $\H$
  corresponding to the parabolic subgroup $W_J$ of $W_r$ is denoted by
  $\H_J$. The {\bf induction functor} from $\H_J$ to $\H_K$ for $J\subseteq
  K\subseteq \delo$ is denoted by $\Ind{J}{K}$. Similarly the
  {\bf restriction functor} from $\H_K$ to $\H_J$ is denoted by $\Res{J}{K}$.
\label{3.4}\end{ddef}
\Medskip

Recall that $\H_J$ is free as an $R$--module with basis $\{T_w\strich w\in
W_J\}$. Moreover, $\H$ is free as a left $\H_J$--module with
basis $\{T_d\strich d\in\D_J\}$, where $\D_J$ denotes the set of
{\bf distinguished right coset representatives} of $W_J$ in $W_r$ which is
the set of all elements of $W_r$ which map the roots in $J$ to
positive roots. Note too, that $\H(\sym{r}) = \H_\Delta$ with $\Delta =
\{\alpha_1,\ldots ,\alpha_r\} = \delo\setminus\{\alpha_0\}$. 

\begin{lemma} Let $\lambda = (\la{1},\la{2})$ be an $a$-bicomposition,
  and let 
$$
A = \{\alpha_0,\alpha_1,\ldots,\alpha_{a-1}\}\cup J_2
$$
where
  $J_2$ is the subset of
  $\{\alpha_{a+1},\alpha_{a+2},\ldots,\alpha_{r-1}\}$ corresponding to
  $\la{2}$. Then $\uap x_\lambda \in\H_A$ and $M^\lambda = 
\Ind{A}{\delo}\uap x_\lambda\H_A$. 
\label{3.5}\end{lemma}

\pf This follows by general arguments since
$\H_A=\H(W_a\times\sym{\lambda^{(2)}})$.
\endproof

We find it convenient to adopt a shorthand for the basis elements of
$M^\lambda$ which appear in Theorem \ref{3.3}.

\begin{notation} Suppose that $\lambda$ is an $a$-bicomposition. Let
$\t = (\ta{1},\ta{2})$ be a $\lambda$--bitableau. We identify the
$\lambda$--bitabloid $\toid{\t}$ and $\uap x_\lambda T_w$,
where $\t^\lambda w$ is the unique row standard $\lambda$--bitableau
which is row equivalent to $\t$.
\label{3.6}\end{notation}
\Medskip

This identification gives us an action of $\H$ 
on the $\lambda$--bitabloids. The discussion after Definition~\ref{2.7} shows 
that this is a $(Q,q)$--analogue of the permutation action of $W_r$ on the
cosets of $W_\lambda$.

Theorem \ref{3.3} says that $M^\lambda$ has a basis which consists of the
distinct $\lambda$--bitabloids. The action of $\H$ on a
$\lambda$--bitabloid $\toid{\t}$ is determined once we know $\{\t\}T_i$
for $i=0,1,\ldots,{r-1}$. From \eqref{3.6} we have the following.

\begin{lemma} Suppose that $0\leq i\leq r-1$ and that both $\t$ and
  $\t s_i$ are row standard (say $\t = \t^\lambda w$). Then 
$$
\{\t\}T_i = \left\{
\begin{alignedat}{2} 
\{\t s_i\}, & &\text{\quad if\quad} \len(ws_i)& = \len(w)+1, \\
q\{\t s_i\}& +(q-1)\{\t\},  &\: \text{\quad if\quad} \len(ws_i)& = \len(w)-1
\text{\quad and\quad} i>0,\\
Q\{\t s_0\}&+(Q-1)\{\t\},  &\: \text{\quad if\quad} \len(ws_i) &= \len(w)-1
\text{\quad and\quad} i=0.  
\end{alignedat} \right.
$$
\label{3.7}\end{lemma}
\Medskip

By considering roots, it is straightforward to determine from the 
tableau~$\t$ which case of the lemma applies.

\begin{lemma} Maintain the notation of Lemma~{\rm\ref{3.7}}. 
Then $\len(ws_i)  = \len(w)+1$ if and only if one of the following holds.
\begin{enumerate}
\item $i,\;i+1\in \t$ and $\row_\t(i) < \row_\t(i+1)$;
\item $-i,\;-i-1\in \t$ and $\row_\t(-i-1) < \row_\t(-i)$;
\item $-i,\; i+1\in \t$; or,
\item $i=0$ and $1\in \ta{2}.$
\end{enumerate}
\label{3.8}\end{lemma}

The cases where $\t s_i$ is not row standard are covered by the next
result.

\begin{lemma} Suppose that $\t = (\ta{1},\ta{2})$ is row standard 
and that $\t s_i$ is not. Then either
\begin{enumerate}
\item $i=0$ and $1\in \ta{1}$ and $\{\t\}T_0 = Q\{\t\}$; or
\item $i>0$ and $i$ and $i+1$ belong to the same row of $\t$ and
  $\{\t\}T_i = q\toid{\t}$; or
\item $i>0$ and $-i$ and $-i-1$ belong to the same row of $\t$ and
  $\toid{\t}T_i = q\toid{\t}$.
\end{enumerate}
\label{3.9}\end{lemma}
\pf It is clear that the only cases where $\t s_i$ is not row standard
are covered by the statement of the lemma.

Consider first the case where $i=0$ and $1\in \ta{1}$. We know that
$\toid{\t} = \uap x_\lambda T_wT_d,$ as in \eqref{3-4}. Moreover, $1$ is
fixed by $d$ because $1\in \ta{1}$. Therefore, $T_0$ commutes with $T_d$, 
and using \eqref{2.3} we get
$$
\toid{\t}T_0 = x_\lambda T_w\uap T_0T_d = Qx_\lambda T_w\uap T_d =
Q\toid{\t}.
$$
This completes the proof of part (i) of the lemma. The proof of the
other parts is similar to that of \cite[3.2(ii)]{DJ1}.
\epf

For future reference, we next collect several results concerning the
action of $\H$ on $M^\lambda$.

\begin{lemma} Suppose that $\t$ is row standard. Assume
  that $i$ is a positive integer in $\ta{2}$. Then
  $\{\t\}T_{i,1}T_0T_{1,i}$ is a linear combination of
  $\lambda$--bitabloids $\toid{\s}$ such that 
  $\s$ is row standard and
\begin{enumerate}
\item $-i\in \s$; and
\item for all integers $j$ with $j\neq \pm i$, we have $j\in \s$ if and
  only if $j\in \t$.
\end{enumerate}
\label{3.10}\end{lemma}
\pf If $i=1$ then the claim follows immediately from Lemmas \ref{3.7}
and~\ref{3.8}. Thus
let $i>1$. Using \eqref{3.7}, \eqref{3.8} and \eqref{3.9}, we see that
$\toid{\t}T_{i-1}$ is a linear combination of one or two
$\lambda$--bitabloids $\toid{\s}$ where $\s$ is row standard and
\begin{enumerate}
\item $i-1\in {\s}^{(2)}$.
\item $i\in \s\Leftrightarrow i-1\in \t$.
\item $j\in \s \Leftrightarrow j\in \t$ for $|j| \neq i-1,i$.
\end{enumerate}
Hence, by induction we may assume that 
$$
\toid{\t}T_{i,1}T_0T_{1,i-1} = \toid{\t}T_{i-1}T_{i-1,1}T_0T_{1,i-1}
$$  
is a linear combination of
$\lambda$--bitabloids $\toid{\s}$ such that  $\s$ is row standard and
\begin{enumerate}
\item $-(i-1) \in {\s}^{(2)}$.
\item $i\in \s \Leftrightarrow i-1\in \t$.
\item $j\in \s\Leftrightarrow j\in \t$ for $|j|\neq i-1,i$.
\end{enumerate}
Therefore, $\toid{\t}T_{i,1}T_0T_{1,i} =
\toid{\t}T_{i-1}T_{i-1,1}T_0T_{1,i-1}T_{i-1}$ is a linear combination
of $\lambda$--bitabloids $\toid{\s}$ as in the statement of the lemma.
\epf

\begin{ddef} Let $M^\lambda_-$ denote the $R$--submodule of $M^\lambda$
  spanned by those $\lambda$--bitabloids $\toid{\t}$ where
$\t=(\t^{(1)},\t^{(2)})$ and $\t^{(2)}$ contains a negative integer.
\label{3.11}\end{ddef} 
\Medskip

Note that Lemmas \ref{3.7} and \ref{3.9} show that
\begin{res}
$M^\lambda_-$ is an $\H(\sym{r})$--submodule of $M^\lambda$.
\label{3.12}\end{res}
\Medskip

Lemma \ref{3.10} now has the following corollary.
\begin{cor} Suppose that $\toid{\t}$ is a $\lambda$--bitabloid. Assume
  that $1,2,\ldots,b \in \ta{2}$. Then
$$
\toid{\t}u_b^- \equiv Q^bq^{b(b-1)/2}\toid{\t}
\pmod{M^\lambda_-}.
$$
\label{3.13}\end{cor}
\pf Note that $Q^bq^{b(b-1)/2}$ is the coefficient of the
identity in $u_b^-$. The corollary follows from multiplying out
$u_b^-$ and applying Lemma \ref{3.10}.
\epf
\begin{lemma} Assume that $\toid{\t}$ is a $\lambda$--bitabloid and that
  all the entries in $\t$ are positive. Suppose that $1\leq i < j\leq
  r$. Then for some integer $k$
$$
\toid{\t}T_{j,i} = q^k\toid{\t s_{j,i}}+m,
$$
where $m$ is a linear combination of $\lambda$--bitabloids $\toid{\s}$
such that $\s$ is row standard and
$\row_{\s}(i) > \row_{\t s_{j,i}}(i) = \row_\t(j)$.
\label{3.14}\end{lemma}
\pf Note first that $js_{j,i} = i$, so $\row_{\t s_{j,i}}(i) =
\row_\t(j)$. Now,
$$
\toid{\t}T_{j-1} = q^\epsilon \toid{\t s_{j-1}} + m_1
$$
where $\epsilon = 0$ or $1$, and either $m_1 = 0$ or
$\row_\t(j-1)>\row_\t(j)$ and $m_1 = (q-1)\toid{\t}$. This proves the lemma
in the case where $i=j-1$. Assume, therefore, that $i<j-1$. 

By induction
$$
\toid{\t s_{j-1}}T_{j-1,i} = q^{k_2}\toid{\t s_{j,i}} + m_2
$$
for some integer $k_2$, where $m_2$ is a linear combination of
$\lambda$--bitabloids $\toid{\s}$ such that $\s$ is row standard and
$\row_{\s}(i) > \row_{\t s_{j-1,i}}(i) =\row_\t(j)$.
Furthermore, if $m_1\neq 0$ then $m_1 = (q-1)\toid{\t}$ and induction
shows that $\toid{\t}T_{j-1,i}$ is a linear combination of
$\lambda$--bitabloids $\toid{\s}$ where $\s$ is row standard and
$$
\row_{\s}(i) \geq \row_{\t s_{j-1,i}}(i)=\row_\t(j-1)>\row_\t(j).
$$ 

Therefore,
$$
\toid{\t}T_{j,i} = \toid{\t}T_{j-1}T_{j-1,i} = q^k\toid{\t s_{j,i}}+m,
$$
for some integer $k$, where $m$ is a linear combination of
$\lambda$--bitabloids $\toid{\s}$ where $\s$ is row standard and
$
\row_{\s}(i) > \row_{\t s_{j,i}}(i).
$
\epf

\begin{cor} Assume that $\toid{\t}$ is a $\lambda$--bitabloid and that
  all the entries in $\t$ are positive. Suppose that $0\leq a\leq
  r$. Then for some integer $k$,
$$
\toid{\t}h_{a,r-a} = q^k\toid{\t w_{a,r-a}}+m,
$$
where $m$ is a linear combination of $\lambda$--bitabloids $\toid{\s}$
($\s$ row standard) such that for all $i$ with $1\leq i\leq r-a$,
$$
\row_{\s}(i) \geq \row_\t(a+i)
$$
with the inequality being strict for at least one $i$.
\label{3.15}\end{cor}

\pf We claim that the following holds. For $0\leq j\leq r-a$
$$ \toid{\t}T_{a+1,1}T_{a+2,2}\cdots T_{a+j,j} 
           = q^{k_1}\toid{\t s_{a+1,1}\cdots s_{a+j,j}}+m
\Number{3-5}$$
where $k_1$ is an integer, and $m$ is a linear combination of
$\lambda$--bitabloids $\toid{\s}$ where $\s$ is row standard and, for 
$1\leq i\leq j$,
$
\row_{\s}(i) \geq \row_\t(a+i),
$
with strict inequality for some $i$ with $1\leq i\leq j$, and
$$
\row_{\s}(i) = \row_\t(i)\quad\text{for}\quad i>a+j.
$$
Equation \eqref{3-5} certainly holds if $j=0$. Since $s_{a+j+1,j+1}$ permutes
only numbers $i$ with $j+1\leq i\leq a+j$, Lemma \ref{3.14} provides the
induction step, so \eqref{3-5} holds.

Since $h_{a,r-a} = T_{a+1,1}T_{a+2,2}\cdots T_{r,r-a}$, by \eqref{2.1},
the corollary is the special case of \eqref{3-5} when $j=r-a$.
\epf



\section{The $(Q,q)$-Schur algebra}

The $q$--Schur algebra \cite{DJ4,DJ5} associated with the Hecke algebra
$\H(\sym r)$ is defined as the endomorphism ring
$$\End_{\H(\sym r)}\(\bigoplus_{\lambda\in\Lambda(n,r)}M^\lambda\)$$
where $\Lambda(n,r)$ is the set of compositions of $r$ into $n$ parts and
$M^\lambda$ is the $\H(\sym r)$--module induced from the trivial module of
the parabolic subalgebra $\H(\sym \lambda)$. Because each partition
in $\Lambda(n,r)$ corresponds to a parabolic subgroup of $\sym r$, a natural 
generalization to a Schur algebra of type $\bf B$ is the algebra
$$ \tilde{\Sch} = \End_\H\(\bigoplus_{J\subseteq\delo}M^J\),\Number{4-1}$$
where $M^J = \Ind{J}{\delo} x_J R$ and
$x_J = \sum_{w\in W_J}T_w$. The algebra $\tilde\Sch$ has been investigated
in \cite{DPS1,G}.

We wish to consider a larger algebra, which contains $\tilde\Sch$ as a
subalgebra, in which we take endomorphisms of induced modules which
correspond to arbitrary reflection subgroups of $W_r$, rather than just the 
parabolic subgroups.

\begin{ddef} The {\bf \Schuh} is the endomorphism ring
$$\Sch_R(n,r) = \End_\H\(\bigoplus_{\lambda\in\bicomp{n,r}}M^\lambda\).$$
\end{ddef}
\Medskip

If $n\geq r$ then all bipartitions of $r$ belong to $\bicomp{n,r}$,
and we see from \eqref{3.2} that the next result holds.

\begin{lemma} If $n\geq r$ then the $(Q,q)$-Schur algebra is Morita 
equivalent to
  $\End_\H\(\bigoplus_\lambda M^\lambda\)$ where the direct sum is now
  over all bipartitions $\lambda$ of $r$.
\label{Morita}\end{lemma}
\Medskip

If $n\leq r$ the $(Q,q)$-Schur algebra $\Sch_R(n,r)$ is Morita equivalent 
to an algebra of the form $e\Sch_R(r,r)e$ for some idempotent $e$ 
in $\Sch_R(r,r)$. Similarly, the algebra $\tilde{\Sch}$ defined in
\eqref{4-1} is Morita equivalent to an algebra of the form $e\Sch_R(n,r)e$.

Henceforth, we fix $n$ and $r$ and let $\Sch = \Sch_R(n,r)$. 

\begin{ddef}
For subsets $K,L,M$ of $\delo$ we set 
$$
\D_{K,L}^M = \D_{K}\cap\D_L^{-1}\cap W_M.
$$
Thus, if $K\subseteq M$ and $L\subseteq M$ then $\D_{K,L}^M$ is the set of 
{\bf distinguished double $(W_K,W_L)$--coset representatives} in $W_M$. 
Similarly, we define $\D_K^M=\D_K\cap W_M$. If $M = \delo$ we usually omit 
the superscript $\delo$. 
\end{ddef}
\Medskip

The next lemma summarizes the properties of the distinguished (double)
coset representatives which we need; the proofs of these results can be
found in~\cite[\S2.7]{Ca}.

\begin{lemma} Suppose that $K,L\subseteq\delo$ and let $d\in\D_{K,L}$. 
\begin{enumerate}
\item Every element of $W_K d W_L$ is uniquely expressible in the form
$xdw$ where $x^{-1}\in\D^K_{K\cap Ld^{-1}}$ and $w\in W_L$; moreover, 
$\len(xdw)=\len(x)+\len(d)+\len(w)$.
\item $W_{Kd\cap L}=d^{-1}W_Kd\cap W_L$; consequently, 
$T_d^{-1}\H_KT_d \cap \H_L=\H_{Kd\cap L}$.            
\item If $K\subseteq L$ then $\D_K=\D^L_K\D_L$; in particular,
$\D_L\subseteq\D_K$.
\end{enumerate}\label{double cosets}
\end{lemma}

As in \cite[3.4]{DJ1} (compare \cite[1.4]{DJ5}) we have the following 
theorem.

\begin{thm}\label{smallschurbasis} The algebra $\tilde{\Sch}$ is free
  as an $R$--module with basis
$$
\{\psi_{I,J}^d\strich I,J\subseteq\delo,\, d\in\D_{I,J}\},
$$
where $\psi_{I,J}^d:x_J\H\to x_I\H$ is the homomorphism given by
$$  \psi_{I,J}^d(x_Jh) = \sum_{w\in W_IdW_J}T_w h,\qquad (h\in\H) .$$
\end{thm}  
\Medskip


We now want to exhibit a generic basis of $\Sch$, that is a basis which is
independent of the choice of the ring $R$ and the values of the
parameters $Q$ and $q$. By construction it is enough to show that 
$$ \Hlm = \Hom_\H(M^\lambda,M^\mu) \Number{4-3} $$
for $\lambda,\mu\in\bicomp{n,r}$ is free as an $R$--module and has a
generic basis. 

First, from the defining relations \eqref{2-1}
for $W_r$ we see that for any $w\in W_r$ the number of factors
equal to $s_0$ is the same in every reduced expression of~$w$. Thus
the following definition makes sense.

\begin{ddef} Suppose that $w\in W_r$ and $\len(w)=a+b$ where a reduced
expression for $w$ has exactly $b$ factors equal to $s_0$. We set 
$$\qlt{w} = q^aQ^b\quad\text{and}\quad\qmlt{w} = q^{-a}Q^{-b}.$$
\end{ddef}
\Medskip

Next we state Frobenius reciprocity and a result which shows that, as
for Hecke algebras of type $\bf A$, induction from parabolic
subalgebras is not only a left adjoint functor to restriction, but a
right adjoint functor as well.


\begin{thm}\label{frob} Let $J\subseteq \delo$, 
and let $M$ be an $\H_J$--module,
  $N$ an $\H$--module. Then the following hold.
\begin{enumerate}
\item $\Hom_\H\(\Ind{J}{\delo}M\,,\,N\) 
\cong \Hom_{\H_J}\(M\,,\,\Res{J}{\delo}N\)$
  where an isomorphism is given by
  restricting a map in $\Hom_\H\(\Ind{J}{\delo}M\,,\,N\)$ to the
  $\H_J$-subspace $M\otimes1$ of 
  \mbox{$\Ind{J}{\delo}M=M\otimes_{\H_J}\H$}. 
\item $\Hom_\H\(N\,,\,\Ind{J}{\delo}M\) \cong
  \Hom_{\H_J}\(\Res{J}{\delo}N\,,\,M\)$ where an isomorphism given by
  sending $\phi\in\Hom_{\H_J}\(\Res{J}{\delo}N\,,\,M\)$ to the map 
  $\hat{\phi}\in\Hom_\H\(N\,,\,\Ind{J}{\delo}M\)$ where $\hat{\phi}$ is 
  defined on $x\in N$ by 
$$
\hat{\phi}(x) = \sum_{d\in\D_J}\qmlt{d} \phi(xT_d^*)\otimes T_d.
$$
\end{enumerate}
\end{thm} 
\pf Part (i) is Frobenius reciprocity which holds for arbitrary
rings and subrings. For part (ii) we refer to \cite[2.6]{DJ1}; the
proof there can be adjusted easily.
\epf

We remark that the inverse to the map
$\phi\mapsto\hat{\phi}$ in part (ii) of the theorem is given by
\begin{equation*}
\theta\mapsto \theta_1\in \Hom_{\H_J}\(\Res{J}{\delo}N\,,\,M\)
\end{equation*}
for $\theta\in \Hom_\H(N\,,\,\Ind{J}{\delo}M)$, where $\theta_1$ is 
defined as follows.
For $x\in N$ we have $x_d\in M$ uniquely defined by
\begin{equation*}
\theta(x) = \sum_{d\in\D_J}x_d\otimes T_d \in \Ind{J}{\delo}M.
\end{equation*}
The map $\theta_d: N \rightarrow M: x\mapsto x_d$ is easily seen to be
$R$--linear and, for $d=1$, even $\H_J$--linear, (in general, $\theta_d$ 
is $T_d^{-1}\H_J T_d$--linear). Then $\theta_1$ is the desired
map.

We shall also use the Mackey Decomposition Theorem \cite[2.7]{DJ1}.

\begin{thm}\label{mackey} Let $I,J$ be subsets of $\delo$ and let $M$ be an
  $\H_J$--module. Then 
$$
\Res{I}{\delo}\Ind{J}{\delo}M = \bigoplus_{d\in\D_{J,I}} 
\Ind{Jd\cap I}{I}\Res{Jd\cap I}{Jd}M\otimes T_d,
$$
considering $M\otimes T_d$ as an $\H_{Jd}$--module, where
$\H_{Jd} = T_d^{-1}\H_JT_d$.
\end{thm}
\Medskip

Now this result taken in conjunction with
Theorem \ref{frob} produces the following intertwining number theorem
(compare \cite[2.8]{DJ1}).
\begin{cor}\label{hommackey} Let $I,J$ and $M$ be as in 
Theorem~{\rm\ref{mackey}} 
and let $N$ be an $\H_I$--module. Let 
$ \frak{H} = \Hom_\H\(\Ind{J}{\delo}M\, ,\, \Ind{I}{\delo}N\)$. Then
\begin{enumerate}\everymath{\displaystyle}
\item
$\begin{array}[t]{ll}
\frak{H} & \cong \bigoplus_{d\in\D_{I,J}}
\Hom_{\H_J}\(M\, ,\, \Ind{Id\cap J}{J}\Res{Id\cap J}{Id}N\otimes T_d\)\\
&\cong \bigoplus_{d\in\D_{I,J}}
\Hom_{\H_{Id\cap J}}\(\Res{Id\cap J}{J}M\, ,\, \Res{Id\cap
  J}{Id}N\otimes T_d\),
\end{array}$
\item $\begin{array}[t]{ll} \frak{H} & \cong \bigoplus_{d\in\D_{J,I}}
  \Hom_{\H_I}\(\Ind{Jd\cap I}{I}\Res{Jd\cap I}{Jd}M\otimes T_d\, ,\, N\)\\
  & \cong  \bigoplus_{d\in\D_{J,I}}
  \Hom_{\H_{Jd\cap I}}\(\Res{Jd\cap I}{Jd}M\otimes T_d\, ,\, \Res{Jd\cap I}I N\).
\end{array}$
\end{enumerate}\end{cor}
\Medskip

\begin{rem}\label{potenzerkl} {\sl We shall apply Corollary 
{\rm \ref{hommackey}} in the special
  case where the modules $M$ and $N$ are the trivial modules 
  $x_J\H_J$ and $x_I\H_I$ respectively. Note that the restriction of
  the trivial module to a parabolic subalgebra is again the
  trivial module. Moreover, the identity map generates the
  endomorphism ring of the trivial module. Corollary {\rm \ref{hommackey}}
  provides ways to exhibit a basis of $\frak{H}$, by choosing the
  identity maps in each summand. Using Theorem {\rm \ref{frob}} the basis 
  elements in parts (i) and (ii) which are labelled by
  $d\in\D_{I,J}$ and $d^{-1}\in\D_{J,I}$ differ by a factor of $\qmlt{d}$.}
\end{rem}
\Medskip

It is an immediate consequence of the way that $W_r$ acts on $V$ that the
next result holds.

\begin{res} Suppose that $w\in W_r$ and $J\subseteq\delo$. Then the following 
are equivalent.
\begin{enumerate}
\item $\alpha_0\in Jw$.
\item $\alpha_0\in J$ and $\alpha_0w=\alpha_0$.
\item $\alpha_0\in J$ and $s_0w=ws_0$ and $\len(s_0w)=\len(w)+1$.
\end{enumerate}\label{new 4.2}\end{res}

>From now on we fix the following notation.
\begin{notation}\label{situation}\sl
Let $a$ and $b$ be integers with $0\le a,b\le r$ and fix an
$a$--bicomposition $\lambda=(\lambda^{(1)}, \lambda^{(2)})$ and a
$b$--bicomposition $\mu=(\mu^{(1)}, \mu^{(2)})$ with
$\lambda,\mu\in\bicomp{n,r}$.
\begin{enumerate}
\item Let $J_1$ be the subset of $\{\alpha_1,\ldots,\alpha_{a-1}\}$
corresponding to $\lambda^{(1)}$ and $J_2$ be the subset of
$\{\alpha_{a+1},\ldots,\alpha_{r-1}\}$ corresponding to $\lambda^{(2)}$.
Similarly, define subsets $I_1$ and $I_2$ corresponding to $\mu^{(1)}$ and 
$\mu^{(2)}$ respectively. Let $J=J_1\cup J_2$ and $I=I_1\cup I_2$.
\item Let $\tA=\delo\setminus\{\alpha_a\}$ if $a\ne r$, and 
$\tA=\delo$ if $a=r$. Similarly, $\tB=\delo\setminus\{\alpha_b\}$ is $b\ne r$,
and $\tB=\delo$ otherwise. Let $A=\tA\setminus\{\alpha_0\}$ and 
$B=\tB\setminus\{\alpha_0\}$.
\end{enumerate}\end{notation}


Note that 
$$ x_J = \sum_{w\in W_J}T_w = x_{J_1}x_{J_2} = x_{J_2}x_{J_1}=x_\lambda. $$
We also remark that 
$J\subseteq A\subseteq \tA\subseteq\delo$ and that
$I\subseteq B\subseteq \tB\subseteq\delo$.

Note that we write $\dou$ instead of $\dou^\delo$.

\begin{lemma}\label{schiebt0}  Let $d\in\dou$, and suppose that $\alpha_0 \in
  \tBdA$. Then $a\geq 1$, $b\geq 1$ and $T_dT_0 =
  T_0T_d$. Moreover for $h\in\H_{\tA}$ we have
\begin{align*}
\uap x_JhT_0 & = Q\uap x_Jh\\
\intertext{and for $h\in \H_{\tB}T_d$ we have}
\ubp x_IhT_0 & = Q\ubp x_Ih.
\end{align*}
\end{lemma}

\pf Since $\alpha_0\in \tB d\cap\tA$, we have $\alpha_0\in\tB\cap\tA$
and $s_0d=ds_0$ by \eqref{new 4.2}. Therefore, $a\ge 1$ and $b\ge 1$ by
\eqref{situation}. Further, since $d\in\D_{\tB,\tA}$, we have
$T_dT_0=T_{ds_0}=T_{s_0d}=T_0T_d$.

Let $h\in\H_{\tA}$. Now  $\tA =                             
\{\alpha_0,\ldots,\alpha_{a-1}\}\cup\{\alpha_{a+1},\ldots\alpha_{r-1}\}$
and hence we may write $h=h_1h_2=h_2h_1$ with 
$h_1\in \H_{\{\alpha_0,\ldots,\alpha_{a-1}\}}$ and
$h_2\in\H_{\{\alpha_{a+1},\ldots,\alpha_{r-1}\}}$.
By \eqref{2.3}, 
$$\uap x_J h T_0=x_J\uap h_1h_2T_0
                =x_Jh_1 \uap T_0h_2
                =Qx_J h_1\uap h_2=Q\uap x_Jh.$$
Similarly, $\ubp x_Ih'T_0=Q\ubp x_Ih'$ for any $h'\in\H_{\tB}$. Since $T_d$
commutes with $T_0$ we have
$\ubp x_I hT_0=Q\ubp x_Ih$ for any $h\in\H_{\tB}T_d$.
\endproof

\begin{ddef} We say that a triple $(d,v,u)$ is {\bf admissible} for
$(\lambda,\mu)$ if $d\in\D_{\tB,\tA}$ and $v\in \D^A_{J,Bd\cap A}$
and $u\in\D^{B}_{I,B\cap Jvd^{-1}}$.
\end{ddef}

Our current aim is to show that $\Hlm$ has a generic basis indexed by the
set of admissible triples for $(\lambda,\mu)$. We first prove three
technical lemmas about admissible triples; after that, the next three lemmas
introduce in turn the elements $d,v$ and $u$ of an admissible triple and
show how they determine a basis of the $R$--module $\Hlm$.

\begin{lemma}\label{intersection1}
Let $(d,v,u)$ be an admissible triple. Then
$B\cap Jvd^{-1}\subseteq Ad^{-1}$.
\end{lemma}

\proof  Let $\alpha\in B\cap Jvd^{-1}$. Then 
$d^{-1}s_\alpha d=v^{-1}s_\beta v$ for
some $\beta\in J$. Now, $v^{-1}s_\beta v\in W_A$ so
$$\len(d)+\len(v^{-1}s_\beta v)=\len(dv^{-1}s_\beta v)
                     = \len(s_\alpha d)
                     =\len(s_\alpha)+\len(d)=1+\len(d).$$
Therefore, $\len(v^{-1}s_\beta v)=1$ and so $\beta v\in A$;
consequently, $\alpha =\beta vd^{-1}\in Ad^{-1}$ proving the lemma.
\endproof

Ultimately we are interested in bases which are indexed not by admissible
triples but by elements of $\D_{I,J}$; this is the point of part (iv) below.

\begin{lemma}Suppose that $(d,v,u)$ is an admissible triple. Then
\begin{enumerate}\item $\len(udv^{-1})=\len(u)+\len(d)+\len(v^{-1})$.
\item $ud\in\D_{I,A\cap Jv}$.
\item $dv^{-1}\in\D_{B,J}$.
\item $udv^{-1}\in\D_{I,J}$.
\end{enumerate}
\label{admissible}\end{lemma}

\proof Note that $d\in\D_{\tB,\tA}\subseteq\D_{B,A}$ and that
$u\in W_B$ and $v^{-1}\in\D^A_{Bd\cap A}$; hence (i) follows from
Lemma~\ref{double cosets}(i).

To prove (ii),  first note that $ud\in \D_I^B\D_B=\D_I$
by Lemma~\ref{double cosets}(iii) because $u\in\D_I^B$ and 
$d\in\D_{\tB,\tA}\subseteq\D_\tB\subset\D_B$. On the other hand, by
Lemma~\ref{intersection1},  
$u\in(\D^B_{B\cap(A\cap Jv)d^{-1}})^{-1}$ and $d\in\D_{B,A\cap Jv}$ so,
by taking $K=B$ and $L=A\cap Jv$ in Lemma~\ref{double cosets}(i), 
$$\len(udw)=\len(u)+\len(d)+\len(w)=\len(ud)+\len(w)$$
for any $w\in W_{A\cap Jv}$.  Therefore, $ud\in(\D_{A\cap Jv})^{-1}$ proving 
(ii). A similar argument proves (iii).

Finally, we prove (iv). By (iii), $dv^{-1}\in\D_B$ so
$udv^{-1}\in\D^B_I\D_B=\D_I$ by Lemma~\ref{double cosets}(iii). 
Furthermore, if $w\in W_J$ then
$$\len(udv^{-1}w)=\len(u)+\len(dv^{-1})+\len(w)
           =\len(udv^{-1})+\len(w)$$
by Lemma \ref{double cosets}(i), since 
$u^{-1}\in\D^B_{B\cap Jvd^{-1}}$
and $dv^{-1}\in\D_{B,J}$ by (iii). Therefore, $udv^{-1}\in\D_J^{-1}$ and
the proof is complete.
\endproof

\begin{lemma} Let $(d,v,u)$ be an admissible triple. Then
$I\cap Jvd^{-1}u^{-1}\subseteq Bu^{-1}$.
\label{intersections}\end{lemma}

\proof Suppose that $\alpha\in I\cap Jvd^{-1}u^{-1}$. Then
$u^{-1}s_\alpha u=dv^{-1}s_\beta vd^{-1}$ for some $\beta\in J$.
Since $u^{-1}s_\alpha u\in W_B$ and $dv^{-1}\in\D_{B,J}$ by
Lemma~\ref{admissible}(iii), 
$$\len(vd^{-1})+\len(u^{-1}s_\alpha u)
      =\len(vd^{-1}u^{-1}s_\alpha u)
      =\len(s_\beta vd^{-1})
      =1+\len(vd^{-1}).$$
Therefore, $\alpha u\in B$ and so the lemma is proved.
\endproof

We are now ready to investigate the Hom--space $\Hlm=\Hom_\H(M^\lambda,M^\mu)$.

\begin{lemma} Suppose that $\lambda$ and $\mu$ are bicompositions in
$\bicomp{n,r}$. Then
$$\Hlm\cong 
     \bigoplus_{d\in\dou}\Hom_{\H_{\tBdA}}\(\Res{\tBdA}{\tA}
       \uap x_J\H_\tA\,,\, \Res{\tBdA}{\tB d}\ubp x_I\H_\tB T_d\).$$
\label{4-8}\end{lemma}

\pf Note that $\{\alpha_0,\ldots,\alpha_{a-1}\}\cup J_2\subseteq\tA$.
Hence, by  Lemma~\ref{3.5} and transitivity of induction,
$M^\lambda\cong \Ind{\tA}{\delo}\uap x_J\H_\tA$. Similarly, we have
$M^\mu\cong\Ind\tB\delo\ubp x_I\H_\tB$.
This is the key observation which allows us to apply Theorems~\ref{frob} 
and~\ref{mackey} to $\Hlm$.
\begin{equation*}\allowdisplaybreaks
\begin{split}
\Hlm & = \Hom_\H(M^\lambda,M^\mu)\\
     & \cong \Hom_\H\(\Ind{\tA}{\delo}\uap x_J\H_\tA\,,\,\Ind{\tB}{\delo}\ubp
     x_I\H_\tB\)\\
     & \cong \Hom_{\H_\tA}\(\uap x_J\H_\tA\,,\,
     \Res{\tA}{\delo}\Ind{\tB}{\delo}\ubp x_I\H_\tB\)\\
     & \cong \bigoplus_{d\in\dou}\Hom_{\H_\tA}\(\uap x_J\H_\tA\,,\,
     \Ind{\tBdA}{\tA}\Res{\tBdA}{\tBd}\ubp x_I\H_\tB T_d\)\\
     & \cong \bigoplus_{d\in\dou}\Hom_{\H_{\tBdA}}\(\Res{\tBdA}{\tA}
     \uap x_J\H_\tA\,,\, \Res{\tBdA}{\tBd}\ubp x_I\H_\tB T_d\)
\end{split}
\end{equation*}
\endproof
 
Recall that $\Delta=\delo\setminus\{\alpha_0\}$.

\begin{lemma} Suppose that $d\in\D_{\tB,\tA}$ and let
$$ \dHlm=\Hom_{\H_\tBdA}\(\Res\tBdA\tA\uap x_J\H_\tA\,,\, 
                         \Res{\tBdA}{\tBd}\ubp x_I\H_\tB T_d\).$$
Then
$$ \dHlm   \cong \bigoplus_{v\in\D_{J,\BdA}^A}\Hom_{\H_{S(v)}}
        \(\Res{S(v)}{Jv}x_JT_vR\,,\,\Res{S(v)}{Bd}(\Ind{I}{B}x_IR)T_d\),$$
where $S(v) = Bd\cap A\cap Jv\subseteq\Delta$. 
\label{4-12}\end{lemma}

\proof
If $\alpha_0\in\tBdA$ then $T_0$ acts on 
the modules $\uap x_J\H_\tA$ and $\ubp x_I\H_\tB T_d$ as 
multiplication by $Q$, by Lemma~\ref{schiebt0}; so,
every $\H_\BdA$--linear map between these
modules is automatically $\H_\tBdA$--linear. 

If $\alpha_0\notin\tBdA$ then $\tBdA = \BdA\subseteq\Delta$ 
by~\eqref{new 4.2}. Therefore,
$$\dHlm=\Hom_{\H_{\BdA}}\(\Res{\BdA}{\tA}\uap x_J\H_\tA\,,\, 
                   \Res{\BdA}{\tBd}\ubp x_I\H_\tB T_d\).$$

Now, by \cite[3.6]{DJ6}, $\uap x_J\H_\tA  \cong x_J\H_A$ as $\H_A$--modules 
and, similarly, $\ubp x_I\H_\tB T_d \cong x_I\H_B T_d$ as 
$T_d^{-1}\H_B T_d$--modules. Therefore, by Theorems~\ref{mackey}
and~\ref{frob},
\begin{equation*}
\begin{split}\allowdisplaybreaks
\dHlm  &\cong  \Hom_{\H_{\BdA}}\(\Res{\BdA}{A}x_J\H_A\,,\,
                                \Res{\BdA}{Bd}x_I\H_B T_d\)\\
      & = \Hom_{\H_{\BdA}}\(\Res{\BdA}{A}\Ind{J}{A}x_JR\,,\,
\Res{\BdA}{Bd}x_I\H_B T_d\)\\
      & \cong \bigoplus_{v\in\D_{J,\BdA}^A}\Hom_{\H_{\BdA}}
\(\Ind{S(v)}{\BdA}\Res{S(v)}{Jv}x_J T_vR\,,\,\Res{\BdA}{Bd}x_I\H_B T_d\),\\
      & \cong \bigoplus_{v\in\D_{J,\BdA}^A}\Hom_{S(v)}
\(\Res{S(v)}{Jv}x_J T_vR\,,\,\Res{S(v)}{\BdA}\Res{\BdA}{Bd}x_I\H_B T_d\).
\end{split}
\end{equation*}
The lemma now follows from the transitivity of restriction and the
observation that $\Ind I B x_IR\cong x_I\H_B$.
\endproof

\begin{lemma} Suppose that $d\in\D_{\tB,\tA}$ and $v\in\D^A_{J,Bd\cap A}$
and let 
$$\vdHlm=\Hom_{\H_{S(v)}}\(\Res{S(v)}{Jv}x_JT_vR\,,\,
                \Res{S(v)}{Bd}(\Ind{I}{B}x_IR)T_d\).$$
Then
$$\vdHlm\cong \bigoplus_{u\in\D^B_{I,B\cap Jvd^{-1}}}      
\Hom_{\H_{S[u]}}\(\Res{S[u]}{Jv}x_JT_vR\,,\,\Res{S[u]}{Iud}x_IT_{ud}R\),$$
where $S[u]=S(v)\cap Iud$.
\end{lemma}

\proof We observe that
$ S(v)d^{-1} = B\cap Ad^{-1}\cap Jvd^{-1}\subseteq \Delta $
determines the parabolic subgroup 
$ W_B\cap d(W_A\cap v^{-1}W_Jv)d^{-1} $
of $W_r$ by Lemma~\ref{double cosets}(ii) since
$d^{-1}\in\D_{\tilde{A},\tilde{B}}\subseteq \D_{A\cap Jv,B}$. Therefore,
\begin{equation*}
\begin{split}
\vdHlm & = \Hom_{\H_{S(v)}}\(\Res{S(v)}{Jv}x_JT_vR\,,\,
                            \Res{S(v)}{Bd}(\Ind{I}{B}x_IR)T_d\)\\
       & \cong \Hom_{\H_{S(v)}}\(\Res{S(v)}{Jv}x_JT_vR\,,\,
                                (\Res{S(v)d^{-1}}{B}\Ind{I}{B}x_IR)T_d\)\\
       &\cong \bigoplus_u
           \Hom_{\H_{S(v)}}\(\Res{S(v)}{Jv}x_JT_vR\,,\, 
                 (\Ind{S[u]d^{-1}}{S(v)d^{-1}}
                 \Res{S[v]d^{-1}}{Iu}x_IT_uR)T_d\),
\end{split}
\end{equation*}
observing that $S[u]d^{-1} = S(v)d^{-1}\cap Iu$, 
where the last isomorphism follows by the Mackey theorem \ref{mackey},
and $u$ runs through $\D^B_{I,S(v)d^{-1}}$. Now 
$$S(v)d^{-1}=B\cap Ad^{-1}\cap Jvd^{-1}=B\cap Jvd^{-1}$$ by 
Lemma~\ref{intersection1} and $ud\in\D^B_{I,S(v)}$ by 
Lemma~\ref{admissible}(ii). Therefore 
$$
{d^{-1}}W_{S(v)d^{-1}}d = W_{S(v)}.
$$ 
Moreover $x_IT_uT_dR$ is a
$T_{d^{-1}}T_{u^{-1}}\H_IT_uT_d$-module, and we have 
$$
(\Ind{S[u]d^{-1}}{S(v)d^{-1}}\Res{S[v]d^{-1}}{Iu}x_IT_uR)T_d
= \Ind{S[u]}{S(v)}\Res{S[u]}{Iud}(x_IT_{ud}R),
$$
since $T_uT_d=T_{ud}$ by Lemma~\ref{admissible}(i). We have shown
$$ \vdHlm \cong \bigoplus_{u\in\D^B_{I,B\cap Jvd^{-1}}}
           \Hom_{\H_{S(v)}}\(\Res{S(v)}{Jv}x_JT_vR\,,\, 
                 \Ind{S[u]}{S(v)}\Res{S[u]}{Iud}x_IT_{ud}R\).
$$
So, Frobenius 
reciprocity and the transitivity of restriction complete the proof.
\endproof
  
Combining the last three lemmas we have shown that 
$\Hlm\cong\bigoplus\uvdHlm$ where the sum is over all admissible
triples $(d,v,u)$ and
$$\uvdHlm = \Hom_{\H_{S[u]}}
              \(\Res{S[u]}{Jv}x_JT_vR\,,\,\Res{S[u]}{Iud}x_IT_{ud}R\).$$
Now, $\uvdHlm$ is a free $R$--module of rank $1$, since both modules 
involved are one dimensional. As a nonzero element in $\uvdHlm$ we may choose
$$ \phi_{I,J}^{(d,v,u)}:x_JT_v\mapsto x_IT_{ud}. $$
By Theorem~\ref{frob}, in the $R$--module $\dHlm$ of Lemma~\ref{4-12}, this map 
corresponds to the map (compare also \cite[3.4]{DJ1})
$$
\hat{\phi}_{I,J}^{(d,v,u)} \in  \Hom_{\H_{\BdA}}\(\Res{\BdA}{A}x_J\H_A\,,\,
\Res{\BdA}{Bd}x_I\H_B T_d\),
$$
given by (noting that $ud\in\D_{I,S(v)}$),
{\everymath{\displaystyle}
$$
\hat{\phi}_{I,J}^{(d,v,u)}(x_JT_vh) = \sum_{w\in W_IudW_{S(v)}}T_wh
\Number{phi hat}$$}%
where $h\in\H_{\BdA}$. The
$\H_{\BdA}$--direct summands of 
$$
\Res{\BdA}{A}x_J \H_A = \bigoplus_{v'\in \D_{J,\BdA}^A}
\Ind{S(v')}{\BdA}\Res{S(v')}{Jv'}x_JT_{v'}R
\Number{killing}$$
with $v'\ne v$ are taken to zero by $\hat{\phi}_{I,J}^{(d,v,u)}$.

We denote the image of $q^{\len(v)}\hat{\phi}_{I,J}^{(d,v,u)}$ in $\Hlm$ 
given by Lemmas~\ref{4-12} and~\ref{4-8} by
$\lmphi$. (The power of $q$ is inserted here because in Lemma~\ref{4-12} we 
used part (ii) of Corollary~\ref{hommackey}; see Remark~\ref{potenzerkl}.)  We 
have proved the following.

\begin{thm}\label{basishom} Let $\lambda,\mu\in\bicomp{n,r}$. 
Then $\Hlm$ is a free $R$--module with basis
$$ \stBlm=\{\lmphi\strich (d,v,u) \text{ an admissible triple for
                              $(\lambda,\mu)$}\}.$$
\end{thm}
\medskip

\begin{cor}\label{basisschur} The $(Q,q)$-Schur algebra $\Sch = \Sch_R(n,r)$ 
is a free $R$--module with basis
$$\stB= \{\lmphi\strich \lambda,\mu\in\bicomp{n,r}\,,\, 
          \text{and $(d,v,u)$ an admissible triple for $(\lambda,\mu)$}\}$$
which does not depend on the commutative ring $R$ or the choice of
the parameter values $Q$ and $q$.
\end{cor}
\Medskip

We call the basis $\stB$ the {\bf standard basis} of $\Sch$. 

We next describe the $\mu$--bitabloids which do not belong to $M^\mu_{-}$ (see
Definition~\ref{3.11}) and are in the support of the image of
$\{\t^\lambda\}=\uap x_J$ under the homomorphism $\lmphi$ for an admissible 
triple $(d,v,u)$ for $(\lambda,\mu)$.

\begin{thm}\label{genmapsto} Let $(d,v,u)$ be an admissible triple for
$(\lambda,\mu)$ and let $c=udv^{-1}$. Then $c\in\D_{I,J}$ and
$$
\lmphi (\toid{\t^\lambda})  \equiv  \ubp\psi_{I,J}^c(x_J) 
                           \equiv \sum_{w\in W_IcW_J}\ubp T_w 
\pmod{\Mmu_-}.
$$
\end{thm}

\pf We have already seen in Lemma~\ref{admissible} that $c=udv^{-1}$ is an 
element of $\D_{I,J}$.

In \eqref{phi hat} we constructed the $\H_{Bd\cap A}$--linear map 
$\hat{\phi}_{I,J}^{(d,v,u)}$ in $\dHlm$, and determined the image of $x_JT_vh$
for $h\in\H_{\BdA}$ under this map.

We apply Lemma \ref{schiebt0} (compare Lemma~\ref{4-12}) to get an 
$\H_{\tBdA}$--linear map 
$\IJhatphi: \Res{\tBdA}{\tA}\uap x_J\H_\tA \longrightarrow
\Res{\tBdA}{\tB .d}\ubp x_I\H_\tB T_d$ such that
$$\IJhatphi(\uap x_J T_vh)= \sum_{w\in W_IudW_{S(v)}}\ubp T_wh
\Number{IJhatphi}$$
for $h\in\H_{\tBdA}$.

According to Lemma~\ref{4-8} we have to trace up this map using 
Theorem~\ref{frob}(ii), to get a map 
\begin{equation*}
\begin{split}
\lmhatphi  \in& \Hom_{\H_\tA}\(\uap x_J\H_\tA\,,\,
\Ind{\tBdA}{\tA}\Res{\tBdA}{\tBd}\ubp x_I\H_\tB T_d\)\\
 &\cong \Hom_{\H_\tBdA}\(\Res{\tBdA}{\tA}\uap x_J\H_\tA\, ,\, 
                        \Res{\tBdA}{\tB .d}\ubp x_I\H_\tB T_d\).
\end{split}
\end{equation*}
Let $m\in \uap x_J\H_\tA$. Then by Theorem~\ref{frob} we have
\begin{equation*}
\lmhatphi(m) =
\sum_{f\in\D_\tBdA^\tA}\qmlt{f}\IJhatphi(mT_f^*)T_f.
\end{equation*}
By \eqref{new 4.2}, $\tBdA\cap\Delta = \BdA$; so
$\D_\tBdA^\tA\cap W_\Delta = \D_\BdA^A$.
Therefore, we may rewrite the last equation as
$$ \lmhatphi(m) = \sum_{f\in\D_\BdA^A}q^{-\len(f)}\IJhatphi(mT_f^*)T_f 
        + \zubar.
\Number{with z}$$
where 
$$ \zubar = \sum_{{\tilde{f}}\in\D_\tBdA^\tA\setminus W_\Delta}
\qmlt{\tilde{f}}\IJhatphi(mT_{\tilde{f}}^*)T_{\tilde{f}}.$$
We claim that $\zubar\in\Mmu_-$.
Observe that $\tilde{f}\in W_\tA$ and $d\in (\D_\tA)^{-1}$, 
so $\len(d\tilde{f}) = \len(d)+\len(\tilde{f})$. Moreover by~\eqref{2.3},
$ \ubp x_I\H_\tB = \ubp x_I\H_B $; therefore
$$
\qmlt{{\tilde{f}}}\IJhatphi(mT_{\tilde{f}}^*)T_{\tilde{f}} 
\in \ubp x_I\H_\tB T_{d{\tilde{f}}} = \ubp x_I\H_B T_{d{\tilde{f}}}.
$$
Now $\tilde{f} \in \D_\tBdA^\tA$ but $\tilde{f} \notin W_\Delta$; so,
$s_0$ is involved in $\tilde{f}$, and in $d\tilde{f}$
as well, since $\len(d\tilde{f}) = \len(d)+\len(\tilde{f})$. 
We apply \cite[1.4]{DJ1} and Lemmas \ref{3.7} and \ref{3.8} to conclude that 
$\qmlt{\tilde{f}}\IJhatphi(mT_{\tilde{f}}^*)T_{\tilde{f}}\in\Mmu_-$.
Consequently, $\zubar\in\Mmu_-$ and \eqref{with z} becomes
$$ \lmhatphi(m) \equiv\sum_{f\in\D_\BdA^A}q^{-\len(f)}\IJhatphi(mT_f^*)T_f 
          \pmod{\Mmu_-}.
\Number{without z}$$

We next investigate $\IJhatphi(mT_f^*)$ when $m=\uap x_J$.
Fix $f\in\D_\BdA^A$ and note that $ f^{-1}\in \D_{\emptyset,\BdA}^A $.
Let $v'$ be the distinguished double coset representative in
$W_Jf^{-1}W_\BdA$. Then $ v'\in\D_{J,\BdA}^A $
and $f^{-1} = gv'$ for some $g\in \D_{(\BdA){v'}^{-1}\cap J}^J$.
However, by \eqref{killing},
$\IJhatphi(\uap x_JT_f^*) = 0$ unless $v' = v$; therefore we assume that
$f^{-1}=gv$ for some $g\in \D_{(\BdA){v}^{-1}\cap J}^J$. By 
Lemma~\ref{double cosets}(i), we have a disjoint union
$$ W_Jf^{-1}W_\BdA = W_JvW_\BdA
            = \bigcup_{g^{-1}\in\D_{S(v)v^{-1}}^J} gvW_\BdA,$$
recalling that $v\in\D^A_{J,Bd\cap A}$ and $S(v)=Bd\cap A\cap Jv$.
Hence, modulo $\Mmu_-$, we can rewrite ~\eqref{without z} when $m=\uap x_J$ as
$$\begin{array}{ll}
  \lmhatphi(\uap x_J) &\equiv \Sum_{g^{-1}\in\D^J_{S(v)v^{-1}}}
               q^{-\len(gv)}\IJhatphi(\uap x_JT_{gv})T_{gv}^*\\
        &\equiv \Sum_{g^{-1}\in\D^J_{S(v)v^{-1}}}
               q^{-\len(v)}\IJhatphi(\uap x_JT_{v})T_{gv}^*\\
        &\equiv \Sum_{g^{-1}\in\D^J_{S(v)v^{-1}}}
              q^{-\len(v)}\ubp\sum_{w\in W_IudW_{S(v)}}T_wT_{gv}^*
\end{array}$$
where the second equality comes from the fact that $g\in W_J$ and the third
equality from \eqref{IJhatphi}.

Using Theorem \ref{frob} we can lift the homomorphism $q^{\len(v)}\lmhatphi$
to the homomorphism $\lmphi$ in $\Hom_\H(\Mla,\Mmu)$ by extending it to 
the module
$\Mla = \Ind{\tA}{\delo}\uap x_J\H_\tA$. In particular the image of the 
generator $\uap x_J$ remains unchanged. Consequently, in order to complete
the proof it suffices to prove the next lemma.

\begin{lemma} Suppose that $(d,v,u)$ is an admissible triple and let
$c=udv^{-1}$. Then
$$\sum_{x\in W_IcW_J}T_x=\sum_{g^{-1}\in\D^J_{S(v)v^{-1}}}
                \sum_{y\in W_Iud W_{S(v)}}T_yT_{gv}^*.$$
\end{lemma}

\proof Since $c\in\D_{I,J}$, by Lemma \ref{double cosets}(i)
$$\sum_{w\in W_IcW_J}T_w 
   = \sum_{\stackrel{y\in W_J}{x^{-1}\in\D^I_{I\cap Jc^{-1}}}}
                T_xT_cT_y
   = \sum_{x^{-1}\in\D^I_{I\cap Jc^{-1}}}
                T_xT_{ud}T_v^*\sum_{y\in W_J}T_y.$$
Recall that $S(v)=Bd\cap A\cap Jv$. Suppose that 
$f\in W_{S(v)v^{-1}}\subseteq W_J$, and let $e=v^{-1}fv$. Then 
$e\in W_{S(v)}\subseteq W_{Bd\cap A}$. Since 
$v\in\D_{J,Bd\cap A}$ we have
$\len(ev^{-1})=\len(e)+\len(v^{-1})$ and
$\len(v^{-1}f)=\len(v^{-1})+\len(f)$; so,
$$T_eT_v^*=T_eT_{v^{-1}}=T_{v^{-1}}T_f=T_v^*T_f.$$
Therefore,
$$T_v^*\sum_{y\in W_J}T_y
   =\sum_{\stackrel{f\in W_{S(v)v^{-1}}}{g\in\D^J_{S(v)v^{-1}}}}
               T_v^*T_fT_g^*
  =\sum_{\stackrel{e\in W_{S(v)}}{g\in\D^J_{S(v)v^{-1}}}}
               T_eT_v^*T_g^*,$$
and consequently,
$$\begin{array}{lr}
\Sum_{w\in W_IcW_J}T_w 
   &= \Sum_{x^{-1}\in\D^I_{I\cap Jc^{-1}}}
       T_xT_{ud}\sum_{\stackrel{e\in W_{S(v)}}{g\in\D^J_{S(v)v^{-1}}}}
               T_eT_{gv}^*\\
   &= \Sum_{g\in\D^J_{S(v)v^{-1}}}
         \( \sum_{\stackrel{e\in W_{S(v)}}{x^{-1}\in\D^I_{I\cap Jc^{-1}}}}
               T_xT_{ud}T_e\)T_{gv}^*.
\end{array}$$
Now $ud\in\D_{I,S(v)}$ by Lemma \ref{admissible}(ii) so by
Lemma~\ref{double cosets}(i) it suffices to show that
$I\cap S(v)d^{-1}u^{-1}=I\cap Jc^{-1}$. However, 
$$\begin{array}{rl}
 I\cap S(v)d^{-1}u^{-1}
          &=I\cap (B\cap Ad^{-1}\cap Jvd^{-1})u^{-1}\\
          &=I\cap (B\cap Jvd^{-1})u^{-1} \\
          &= I\cap Jvd^{-1}u^{-1},
\end{array}$$
by Lemma~\ref{intersection1} and Lemma~\ref{intersections} respectively.
\endproof

\begin{rem}\label{groupcase}
{\sl In the group algebra case, that is specializing both $Q$
    and $q$ to $1$, one sees easily that 
$$
\uap x_J = \sum_{w\in W_\lambda}w \in RW_r.
$$
In this case $\lmphi$ takes the generator $\uap x_J$ of the module
$\Mla$ to 
$$
\lmphi(\uap x_J) = \sum_{w\in W_\lambda cW_\mu}w,
$$
where $c=udv^{-1}$ and the subgroups $W_\lambda$ and $W_\mu$ are defined in
\eqref{W-lambda}.}
\end{rem}\Medskip

By \eqref{new 4.2}, $\D_{\tB,\tA}\cap W_\Delta = \D_{B,A}^\Delta.$
Therefore, $\Hlm$ in \eqref{4-8} splits into two $R$--modules $\Hlm^+$ and
$\Hlm^-$ where
\begin{equation*}
\Hlm^+  = \bigoplus_{d\in\D_{B,A}^\Delta}
\Hom_{\H_\tBdA}\(\Res{\tBdA}{\tA}\uap x_J\H_\tA\, ,\,
\Res{\tBdA}{\tBd}\ubp x_I\H_\tB T_d\)
\end{equation*}
and
\begin{equation*}
\Hlm^-  = \bigoplus_{d\in\D_{B,A}\setminus W_\Delta}
\Hom_{\H_\tBdA}\(\Res{\tBdA}{\tA}\uap x_J\H_\tA\, ,\,
\Res{\tBdA}{\tBd}\ubp x_I\H_\tB T_d\).
\end{equation*}

We apply Lemma \ref{schiebt0} and Theorem \ref{frob} again to get
\begin{equation*}\begin{split}
\Hlm^+ & \cong \bigoplus_{d\in\D_{B,A}^\Delta}
\Hom_{\H_\BdA}\(\Res{\BdA}{A}\uap x_J\H_A\, ,\,
\Res{\BdA}{B d}\ubp x_I\H_B T_d\)\\
& \cong \bigoplus_{d\in\D_{B,A}^\Delta}
\Hom_{\H_A}\(\uap x_J\H_A\, ,\, \Ind{\BdA}{A}
\Res{\BdA}{Bd}\ubp x_I\H_B T_d\)\\
& \cong \Hom_{\H_\Delta}\(x_J\H_\Delta\, ,\, x_I\H_\Delta\).
\end{split}\end{equation*}

The last $\Hom$-set has basis 
$ \stB_\Delta = \{\IJtphi \strich c\in \D_{I,J}^\Delta\}$,
where the map $\IJtphi$ is defined by
$$\IJtphi(x_Jh)  = \sum_{w\in W_IcW_J}T_wh \qquad(h\in\H_\Delta).$$

We observe that the map $\IJtphi$ involves the same
double coset $W_IcW_J$ as the map $\psi_{I,J}^c$ in Theorem
\ref{genmapsto}. Moreover, it follows that
$$
\stBlm^+ = \{\lmphi\in\stBlm\strich d\in \D^\Delta_{B,A}
              \text{ and $(d,v,u)$ admissible for $(\lambda,\mu)$}\}
$$
is a basis of $\Hlm^+$ since $\stBlm^+$ and $\stB_\Delta$ have the same
cardinality. This proves the next lemma (it is also easily checked directly).

\begin{lemma}\label{blackbasis} Let $\lambda,\mu\in\bicomp{n,r}$ and
let $I$ and $J$ be as in~\eqref{situation}. Then
$$\D_{I,J}^\Delta = 
\{udv^{-1}\strich  (d,v,u) \text{ an admissible triple for $(\lambda,\mu)$
                    and $d\in\D^\Delta_{B,A}$}\}$$
\end{lemma}\Medskip

In view of the last result we make the following
definition.

\begin{ddef}\label{doublecosetmap} 
  Given $c\in\D^\Delta_{I,J}$ let $\lmcphi\in\stBlm$ be the map
  $\lmphi$ where $(d,v,u)$ is the admissible triple such that $c=udv^{-1}$.
\end{ddef}\Medskip

Observe that Theorem~\ref{genmapsto} implies that for an
admissible triple $(d,v,u)$ and $c=udv^{-1}$
$$
\lmphi(\toid{\t^\lambda})\equiv \sum_{b\in\D_{Ic\cap
    J}^J}\toid{\t^\mu}T_{cb}\pmod{\Mmu_-}. 
\Number{tableaux sums}$$
>From this the first part of the next result follows easily. The
second part is a special case of Theorem~\ref{genmapsto}.
 
\begin{cor}\label{blackd}  Let $(d,v,u)$ be admissible, and let 
$c = udv^{-1}$.
\begin{enumerate}
\item  If $c\notin\D^\Delta_{I,J}$ 
then $\lmphi(\toid{\t^\lambda})\equiv 0\pmod{\Mmu_-}$.
\item If $c\in\D^\Delta_{I,J}$, then
$$
 \lmcphi(\toid{\t^\lambda})\equiv \ubp\sum_{w\in W_IcW_J}T_w \pmod{ \Mmu_-}.$$
\end{enumerate}
\end{cor}
\Medskip  

 

\setcounter{equation}{0}
\section{$(Q,q)$--tensor space and $(Q,q)$--Weyl modules}

In view of Lemma \ref{Morita} we assume hereafter that $n\geq r$. We denote
the bicomposition $((-),(1^r))$ by $\om$. Thus $\om\in\bicomp{n,r}$,
since $n\geq r$. 

\begin{ddef}
{\bf $(Q,q)$-tensor space} is the $R$--module
$$ \E = \E_R(n,r) = \bigoplus_{\lambda\in\bicomp{n,r}}\Mla.  $$
\end{ddef}

\begin{ddef}\label{idempdef} Let $\lambda$ be an $a$--bicomposition in 
$\bicomp{n,r}$. 
\begin{enumerate}
\item The homomorphism $\llphi$ is the identity map on $\Mla$ and
  maps $\Mmu$ to zero when $\lambda\ne\mu$.
\item The homomorphism $\lophi$ is premultiplication of $\H =
  M^\om$ by $\uap x_\lambda$ and maps $\Mmu$ to zero when $\mu\ne\omega$.
\end{enumerate}
\end{ddef}\Medskip

Thus $\llphi$ is the projection of $\E$ onto $\Mla$ and $\lophi$ is
the canonical epimorphism of $\H$ onto the cyclic $\H$-module
$\Mla$. In particular $\oophi$ is the identity on $\H$ and maps $\Mla$
to zero for $\lambda \neq \om$. From the proof of Theorem~\ref{genmapsto}, 
using in particular \eqref{with z}, one sees easily that the following holds.

\begin{lemma}\label{llundlo} Let $\lambda$ be a bicomposition. 
\begin{enumerate}
\item The restriction of $\llphi$ to $M^\lambda$ equals 
$\llphi^{(1,1,1)}$.
\item The restriction of $\lophi$ to $\H$ equals $\lophi^{(1,1,1)}$.
\end{enumerate}
\end{lemma}\Medskip

By construction $\{\llphi \strich \lambda\in\bicomp{n,r}\}$ is a set of
orthogonal idempotents of $\Sch$ whose sum is the identity element of
$\Sch$. Observe that
\begin{equation}
\llphi\lophi  = \lophi = \lophi\oophi,
\end{equation} 
hence, as in \cite[2.5]{DJ5}, we see that postmultiplication by $\lophi$ 
embeds the left ideal $\Sch\llphi$ of $\Sch$ in the left ideal $\Sch\oophi$. 
Hence, $\Sch$ acts faithfully on $\Sch\oophi$ and, when $R$ is a field, every 
irreducible left $\Sch$-module occurs as a composition factor of the left 
ideal $\Sch\oophi$ of $\Sch$. We have the following 
(compare \cite[2.10 and 2.6]{DJ5}). 
 
\begin{lemma}\label{tensorspace} The Hecke algebra $\H$ is canonically
  isomorphic to $\oophi\Sch\oophi$ and acts on $\Sch\oophi$ as a set of
  $\Sch$--linear maps. The $\H$--submodule $\lophi\H$ of $\Sch\oophi$ is
  isomorphic to $\Mla$ and $(Q,q)$-tensor space $\E$ is
  isomorphic to the left ideal $\Sch\oophi$ of $\Sch$ as an
  $(\Sch$,$\H$)--bimodule.
\end{lemma}

\pf We have canonical isomorphisms
$$
\H \cong \End_\H(\H) \cong \Hom_\H(x_\om\H\, , \,x_\om\H) \cong
\oophi\Sch\oophi.
$$
Identifying $\H$ and $\oophi\Sch\oophi$ we see   that $\H$
acts on $\Sch\oophi$ on the right as a set of $\Sch$--linear maps. Also,
$\lophi h\mapsto \uap x_\lambda h$ for $h\in\H$
gives an $\H$--isomorphism between $\lophi\H$ and $M^\lambda$.

By premultiplying $\Sch\oophi$ by
$\sum_{\lambda\in\bicomp{n,r}}\llphi$, which is the identity of $\Sch$, we
obtain 
$$\Sch\oophi=\bigoplus_{\lambda\in\bicomp{n,r}}M^\lambda=\E.
\endproof$$

\begin{ddef}\label{weightspacedef}   Let $U$ be a left $\Sch$-module, and let
  $\lambda\in\bicomp{n,r}$. Then the $R$--submodule $U_\lambda = \llphi U$ 
  of $U$ is the {\bf weight space} of $U$ of {\bf weight} $\lambda$.
\end{ddef}\Medskip

Note that $\llphi U$ is free as an $R$--module. We have
\begin{equation}
U = \bigoplus_{\lambda\in\bicomp{n,r}}\llphi U,
\end{equation}
since $\{\llphi\strich \lambda\in\bicomp{n,r}\}$ is a set of
orthogonal idempotents whose sum is the identity of $\Sch$. This decomposition 
of $U$ is called the {\bf weight space decomposition} of $U$, (compare
\cite[2.13]{DJ5}). 

\begin{lemma}\label{gewichtszerl} The weight space decomposition of
  $(Q,q)$--tensor space $\E$ is given as 
$$
\E = \bigoplus_{\lambda\in\bicomp{n,r}}\Mla.
$$
More generally, if $U$ is any $R$--submodule of $\E$ then its weight space
decomposition is
$$
U = \bigoplus_{\lambda\in\bicomp{n,r}}U\cap\Mla.
$$
\end{lemma}

Let $\lambda\in\bicomp{n,r}$ and recall the definition of $z_\lambda$ 
from~\eqref{z-lambda def}. With the above identifications we have
$$
z_\lambda = \lophi h_{a,r-a}T_{\hat{\pi}_\lambda}u_{r-a}^-\hat{y}_\lambda.
\Number{z-lambda}$$

Contrast the next definition with our definition of the Specht module
$S^\lambda = z_\lambda\H$.

\begin{ddef}\label{weylmoduledef} Suppose that $\lambda$ is an 
$a$-bicomposition of $r$.
\begin{enumerate}
\item Let $W^\lambda = \Sch z_\lambda$. 
\item Let $L^\lambda = \Sch u_{r-a}^-\hat{y}_\lambda$.
\end{enumerate}
We call $W^\lambda$ a {\bf $(Q,q)$-Weyl module}.
\end{ddef}\Medskip

As in \cite[3.9]{DJ5} we have the following.

\begin{res}\label{moritaweyl} If $\lambda$ and $\mu$ are associated
  bicompositions then $W^\lambda\cong W^\mu$ and $L^\lambda\cong L^\mu$.
\end{res}\Medskip

We next define a bilinear form $\<\ ,\ \>_\lambda$ on $\Mla$ by specifying 
its values at pairs of basis elements, which are given in \eqref{3.6} by 
$\lambda$--bitabloids. 
For such a pair $\toid{\t},\toid{\t' }$, where $\t = \t^\lambda w$ is 
row standard and $w\in W_r$, we set
\begin{equation}
\<\toid{\t},\toid{\t' }\>_\lambda = 
\begin{cases} \qlt{w} & \text{if $\toid{\t}=\toid{\t' }$}\\
                 0    & \text{if $\toid{\t}\neq\toid{\t' }$}.
\end{cases}\end{equation}

Extend these bilinear forms to a non-degenerate bilinear form $\<\ ,\ \>$
on $(Q,q)$-tensor space $\E = \bigoplus\Mla$ by specifying that
$\Mla$ and $\Mmu$ are orthogonal when $\lambda \neq \mu$.

Recall that the anti--automorphism $^*$ on $\H$ is defined by $T_w^* =
T_{w^{-1}}$ for $w\in W_r$.
As in \cite[4.1]{DJ5} for all $h\in\H$, $x,y\in\E$
\begin{equation}
\<xh,y\> = \<x,yh^*\>,
\end{equation}
and consequently $\E$ is self--dual as an $\H$-module. Here $M^* =
\Hom_\H(M,R)$ is the dual of an $\H$-module $M$, the right action of
$\H$ on $M^*$ being given by $fh(x) = f(xh^*)$,
for $f\in M^*$, $h\in\H$ and $x\in M$. The discussion in \cite[section
1]{DJ5} shows that $*$ extends to an anti--isomorphism of $\Sch$, which we
also denote by~$*$. This allows us to define the dual of a left
$\Sch$--module. We then have the following.

\begin{lemma}\label{bilforminv} Let $x,y\in\E$, $h\in\H$ and
  $s\in\Sch$. Then  
\begin{align}
\<sx,y\> & = \<x,s^*y\>\tag{\rm i}\\
\<xh,y\> & = \<x,yh^*\>\tag{\rm ii}
\end{align}
In particular, $(Q,q)$--tensor space $\E$ is self--dual as a left
$\Sch$--module.
\end{lemma}\Medskip

Since $(u_{r-a}^-\hat{y}_\lambda)^* = u_{r-a}^-\hat{y}_\lambda$ we can
contract the bilinear form which is given by restricting $\<\ ,\ \>$ to
$L^\lambda$ to obtain a bilinear form $\Ll \ ,\ \Gg $ on $L^\lambda$
such that 
\begin{equation*}\begin{split}
\Ll s_1u^-_{r-a}\hat{y}_\lambda\, ,\, s_2u^-_{r-a}\hat{y}_\lambda\Gg & =
\Ll s_1\oophi u^-_{r-a}\hat{y}_\lambda\, ,\, s_2\oophi
u^-_{r-a}\hat{y}_\lambda\Gg \\
& = \<s_1\oophi\, ,\, s_2u^-_{r-a}\hat{y}_\lambda\>,
\end{split}\end{equation*}
for all $s_1,s_2\in\Sch$, (compare \cite[section 5.5]{Gr2} and
\cite[4.2]{DJ5}). The proof of \cite[4.4]{DJ5} now gives the following
theorem.

\begin{thm}[The Submodule Theorem]\label{submodulethm} 
Suppose that $R$ is a field and
  $\lambda \in\bicomp{n,r}$. Let $U$ be an $\Sch$--submodule of
  $L^\lambda$. Then $W^\lambda\subseteq U$ or $U\subseteq W^{\lambda\perp}$.
\end{thm}\Medskip

\begin{cor}\label{weylhead} Suppose that $R$ is a field and
  $\lambda\in\bicomp{n,r}$. Then $W^\lambda\cap W^{\lambda\perp}$ is the
  unique maximal submodule of $W^\lambda$, and the quotient
  $W^\lambda/(W^\lambda\cap W^{\lambda\perp})$ is an absolutely
  irreducible self--dual $\Sch$-module. 
\end{cor}

\pf The corollary will be an immediate consequence of the Submodule Theorem
once we prove that the generator $z_\lambda$
of $W^\lambda$ is anisotropic. We have using \eqref{z-lambda} and
\eqref{3.6} 
\begin{equation*}\begin{split}
\Ll z_\lambda,z_\lambda\Gg  & =
\<\lophi h_{a,r-a}T_{\hat{\pi}_\lambda}\, ,\, 
\lophi h_{a,r-a}T_{\hat{\pi}_\lambda}u_{r-a}^-\hat{y}_\lambda\>\\
& = \<\toid{\t}\, ,\, \toid{\t}u_{r-a}^-\hat{y}_\lambda\>,
\end{split}\end{equation*}
where $\t = \hat{\t}_\lambda = \tlam w_{a,r-a}\hat{\pi}_\lambda$
(see Definitions \ref{2.7} and \ref{2.12}). Hence , by Corollary~\ref{3.13},
$$
\Ll z_\lambda,z_\lambda\Gg  =
Q^{r-a}q^{(r-a)(r-a-1)/2}\<\toid{\t},\toid{\t}\hat{y}_\lambda\>.
$$
But $\toid{\t}\hat{y}_\lambda = \toid{\t} + v$, where $v$ is a linear
combination of $\lambda$--bitabloids distinct from $\toid{\t}$
(cf. \cite[4.1]{DJ1}). Therefore
$$
\Ll z_\lambda,z_\lambda\Gg  =
Q^{r-a}q^{(r-a)(r-a-1)/2}\<\toid{\t},\toid{\t}\> \neq 0.
$$
\epf

\begin{ddef}\label{sirred} Suppose that $R$ is a field and
  $\lambda\in\bicomp{n,r}$. Let $F^\lambda$ be the irreducible
  $\Sch_R(n,r)$-module $W^\lambda/(W^\lambda\cap W^{\lambda\perp})$. 
\end{ddef}\Medskip

\begin{thm}\label{decmatrix} Suppose that $R$ is a field and that
  $\lambda,\mu\in\bicomp{n,r}$. We have $F^\lambda\cong F^\mu$ if and
  only if $\lambda$ and $\mu$ are associated bicompositions. Thus 
$$
\{F^\lambda\strich \lambda \text{ is a bipartition of } r\}
$$
is a set of non--isomorphic absolutely irreducible self--dual
$\Sch_R(n,r)$-modules.

If $\lambda$ and $\mu$ are bipartitions of $r$ then let
$d_{\lambda\mu}$ be the multiplicity of $F^\mu$ as a composition
factor of $W^\lambda$. If the bipartitions of $r$ are ordered
lexicographically then the matrix $(d_{\lambda\mu})$ is upper
unitriangular.
\end{thm}
\pf This theorem is proved in exactly the same way as \cite[4.11 and
4.13]{DJ5} (note \cite[3.7]{DJ6}).
\epf


\setcounter{equation}{0}
\section{The semistandard basis theorem}

Hereafter, we assume that $\lambda,\mu\in\bicomp{n,r}$, with $\lambda$
being an $a$-bipartition.

\begin{ddef}\label{tabloftype}
\begin{enumerate}
\item A {\bf $\lambda$--bitableau of type $\mu$} is an array $\T=
  \btabl$ of integers obtained from the diagram $[\lambda]$ by
  replacing each cross by a non-zero integer according to the
  following restrictions. For $1\leq i\leq n$, the number of entries
  $j$ with $j=i$ is equal to the $i$th part of $\mue{1}$; and the
  number of entries $j$ with $|j|=n+i$ is equal to the $i$th part of
  $\mue{2}$. For $k=1,2$, $\Ta{k}$ denotes the array of integers
  replacing the crosses in $\la{k}$. We denote the set of
  $\lambda$--bitableaux of type $\mu$ by $\Bit$. Note that all entries
  less than $n$ in $\T\in\Bit$ are positive.
\item An element $\T = \btabl$ of $\Bit$ is {\bf positive} if all the
  entries in $\Ta{1}$ are greater than $0$ and all the entries in
  $\Ta{2}$ are greater than $n$. We denote the set of positive
  elements of $\Bit$ by $\Bitplus$.
\item A {\bf semistandard $\lambda$--bitableau of type $\mu$} is an
  element of $\Bitplus$ in which the entries are weakly increasing
  along each row and strictly increasing down each column. We denote
  the set of semistandard $\lambda$--bitableaux of type $\mu$ by $\Bitss$.
\end{enumerate}
\end{ddef}\Medskip

\begin{ex}\label{beisp1} If $r=n=7$ and $\lambda = ((3,2),(1,1))$ and
  $\mu=((3,1),(2,1))$ then
\begin{alignat*}{3}
\text{\ } \qquad \T_1 &&=  \left(\begin{matrix}    1&8&9\\
                                           2&1&\ \  
\end{matrix}\right.
                   \quad &, \quad
                  \left.\begin{matrix}         \  1\\
                                               -8  \end{matrix}
\right ) &&\in\Bit,\\
\text{and }  \qquad \T_2 &&= \left(\begin{matrix} 1&1&1\\
                                           2&8&\     
\end{matrix}\right.
                    \quad &, \quad
                  \left.\begin{matrix}        \ 8\\
                                          \ 9 
                                        \end{matrix}
            \right )&& \in \Bitss.
\end{alignat*}
\end{ex}
\Medskip

Note that there are no positive $\lambda$--bitableaux of type $\mu$ if
$\lambda$ is an $a$-bipartition and $\mu$ is a $b$-bicomposition of
$r$ with $a<b$.

\begin{ddef}\label{ordnung}
Given $\T\in\Bitplus$, let $\T(i)$ equal the integer
which appears in $\T$ in the place which is occupied by $i$ in
$\hat{\t}_\lambda$. Choose a total order $<$ on $\Bitplus$ such that if $\A$
and $\B$ are elements of $\Bitplus$ then  $\A<\B$ if one of the following
holds.
\begin{enumerate}
\item For all $j,k$.
\begin{multline*}
\#\{i\strich \A(i) = j\text{ and $i$ belongs to column $k$ of
  }\hat{\t}_\lambda\}\\
=\#\{i\strich \B(i) = j\text{ and $i$ belongs to column $k$ of
  }\hat{\t}_\lambda\},
\end{multline*}
and $\A(i)<\B(i)$ for the smallest integer $i$ such that $\A(i)\neq \B(i)$.
\item $\Sum_{i=1}^{r-a}\A(i) < \sum_{i=1}^{r-a}\B(i)$.
\item $\Sum_{i=1}^{r-a}\A(i) = \sum_{i=1}^{r-a}\B(i)$, and for all $j,k$
\begin{multline*}
\#\{i\strich \A(i)\leq j\text{ and $i$ belongs to the first $k$ columns of }
\hat{\t}_\lambda\}\\
\geq \#\{i\strich \B(i)\leq j\text{ and $i$ belongs to the 
first $k$ columns of }\hat{\t}_\lambda\}.
\end{multline*}
\end{enumerate}
\end{ddef}\Medskip

Let $V(\lambda,\mu)$ be the free $R$-module spanned by all $\T\in\Bit$.
Given $\T\in\Bitplus$, define $E_\T$ to be the element of $V(\lambda,\mu)$ obtained
by summing over those $\T'$ which are row equivalent to $\T$, and then applying 
the signed column symmetrizer in a fashion similar to that 
in~\cite[8.1.11]{JK}.

\begin{ex}\label{beisp2} If 
$$
\T =  \left(\begin{matrix} 1&1\\ 2   \end{matrix}\  ,\ 
            \begin{matrix} 5\\6 \end{matrix}\ \right ),
$$
then
$$
E_\T = \left(\begin{matrix}  1&1\\ 2&\   \end{matrix}\ ,\ 
              \begin{matrix} 5\\6\\\end{matrix}\ \right)
     - \left(\begin{matrix}  2&1\\ 1&\   \end{matrix}\ ,\ 
              \begin{matrix} 5\\6\\\end{matrix}\ \right)
     - \left(\begin{matrix}  1&1\\ 2&\   \end{matrix}\ ,\ 
              \begin{matrix} 6\\5\\\end{matrix}\ \right)
     + \left(\begin{matrix}  2&1\\ 1&\   \end{matrix}\ ,\ 
              \begin{matrix} 6\\5\\\end{matrix}\ \right)
$$
\end{ex}\Medskip

We shall need the following result, which follows easily from the
corresponding theorem \cite[8.1.16]{JK} for Weyl modules for general linear 
groups.

\begin{thm}\label{dominantbase} The $R$--submodule of $V(\lambda,\mu)$ spanned by
  $\{E_\T\strich \T\in\Bitplus\}$ is a free $R$--module with basis
  $\{E_\T\strich \T\in\Bitss\}$.
\end{thm}\Medskip 

\begin{ddef}\label{sigmamap} \begin{enumerate}
\item If $\T = \btabl \in \Bit$ then let $\ep(\T)=1$ if an even number
  of entries in $\Ta{2}$ are negative , and $\ep(\T) = -1$, otherwise.
\item Let $\sigma$ be the linear map on $V(\lambda,\mu)$ which sends each
  $\T\in\Bit$ to 
$$
\sigma(\T) = \sum\ep(\T')\T,
$$
where the sum ranges over all those
  $\T'\in \Bit$ which are obtained from $\T$ by changing the sign of
  some of the entries in $\T$ which are greater than $n$. 
\end{enumerate}
\end{ddef}\Medskip

\begin{ex}\label{beisp3} If
$$
\T =
\left(\begin{matrix}    1&1\\
                        2&\   
\end{matrix}\right.
                   \,  ,\, 
   \left.\begin{matrix} 5&6\\
                       \ &\ \end{matrix}\right ),
$$
then 
\begin{multline*}
\sigma(\T) = 
\left(\begin{matrix}    1&1\\
                        2&\   
\end{matrix}\right.
                   \,  ,\, 
   \left.\begin{matrix} 5&6\\
                       \ &\ \end{matrix}\right )
-
\left(\begin{matrix}    1&1\\
                        2&\   
\end{matrix}\right.
                   \,  ,\, 
   \left.\begin{matrix} -5&6\\
                       \ &\ \end{matrix}\right )\\
-
\left(\begin{matrix}    1&1\\
                        2&\   
\end{matrix}\right.
                   \,  ,\, 
   \left.\begin{matrix} 5&-6\\
                       \ &\ \end{matrix}\right )+
\left(\begin{matrix}    1&1\\
                        2&\   
\end{matrix}\right.
                   \,  ,\, 
   \left.\begin{matrix} -5&-6\\
                       \ &\ \end{matrix}\right ).
\end{multline*}
\end{ex}\Medskip

Theorem \ref{dominantbase} immediately gives the following corollary.

\begin{cor}\label{domsignbase}The $R$--module generated by
  $\{\sigma(E_\T)\strich \T\in\Bitplus\}$ is free with basis  
$\{\sigma(E_\T)\strich \T\in\Bitss\}$.
\end{cor}\Medskip

To facilitate our calculations, we introduce for fixed bicompositions
$\lambda,\mu\in\bicomp{n,r}$ three  linear
maps, $\alpha$, $\beta$, and $\gamma$ from $\Mmu$ into $V(\lambda,\mu)$.

\begin{ddef}\label{maps} Let $\alpha,\beta,\gamma$ be the linear
  transformations from $\Mmu$ into $V(\lambda,\mu)$ which are given as
  follows. Suppose that $\t= (\ta{1},\ta{2})$ is a 
  $\mu$--bitableau. The $\lambda$--bitableau $\alpha\toid{\t}$ 
  of type $\mu$ is obtained from $\tlam$ and $\t$ as
  follows. For $1\leq i\leq r$
\begin{enumerate}
\item replace the entry $i$ in $\tlam$ by $j$  if $i$
  or $-i$ occurs in row $j$ of $\ta{1}$, and
\item replace the entry $i$ in $\tlam$ by $n+j$ (respectively $-n-j$)
  if $i$ (respectively $-i$) occurs in row $j$ of $\ta{2}$.
\end{enumerate}
The definitions of $\beta\toid{\t}$ and $\gamma\toid{\t}$ are
obtained in a similar way, replacing $\tlam$ by $\hat{\t}^\lambda$ and
$\hat{\t}_\lambda$ respectively. 
\end{ddef}\Medskip
Observe that these maps are independent of the choice of tableau $\t$ in
$\toid{\t}$ and so are well defined.

\begin{ex}\label{beisp4} Assume that $r=n=6$, and let $\lambda =
  ((2^2),(2))$, $\mu = ((3,1),(1^2))$. Then 
\begin{equation*}
\tlam  =
\left(\begin{matrix}    1&2\\
                        3&4   
\end{matrix}\right.
                   \,  ,\, 
   \left.\begin{matrix} 5&6\\
                       \ &\ \end{matrix}\right )\quad
\hat{\t}^\lambda =
\left(\begin{matrix}    3&4\\
                        5&6   
\end{matrix}\right.
                   \,  ,\, 
   \left.\begin{matrix} 1&2\\
                       \ &\ \end{matrix}\right )\quad
\hat{\t}_\lambda =
\left(\begin{matrix}    3&5\\
                        4&6   
\end{matrix}\right.
                   \,  ,\, 
   \left.\begin{matrix} 1&2\\
                       \ &\ \end{matrix}\right ).
\end{equation*}
If 
$$ \t = \left(\begin{matrix} 3&4&6\\ 1&\ \ \end{matrix}\right.
                   \,  ,\, \left.\begin{matrix} -2\\
                       5\end{matrix}\right )
$$ 
then
$$
\alpha\toid{\t} = \left(\begin{matrix} 2&-7\\ 1&1\end{matrix}\right. \, ,\, 
   \left.\begin{matrix} 8&1\\ \ &\ \end{matrix}\right )\quad,\quad
\beta\toid{\t} = \left(\begin{matrix} 1&1\\ 8&1 \end{matrix}\right. \, ,\, 
   \left.\begin{matrix} 2&-7\\ \ &\ \end{matrix}\right )
$$
and
$$
\gamma\toid{\t}= \left(\begin{matrix}    1&8\\ 1&1   \end{matrix}\right.
\,  ,\, \left.\begin{matrix} 2&-7\\ \ &\ \end{matrix}\right ).
$$
\end{ex}\Medskip

Note that all three maps defined above have inverses. Given $\T\in\Bit$
we define the $\mu$--bitabloid $\toid{\t} = \alpha^{-1}\T$ as follows. If
the place occupied by $i$ in $\t^\lambda$ is occupied by $j$
(respectively by $-j$) in $\T$ put $i$ (respectively $-i$) in row $j$
of $\t$, counting the rows of $\ta{1}$ as row $1,2,\ldots$ and the rows
of $\ta{2}$ as row $n+1,n+2,\ldots$, and then take the bitabloid
$\toid{\t}$ containing the bitableau $\t$. For the maps $\beta$ and
$\gamma$ use $\hat{\t}^\lambda$ and $\hat{\t}_\lambda$ respectively 
instead of $\t^\lambda$.

\begin{lemma}\label{inversemaps} The maps $\alpha$, $\beta$ and
  $\gamma$ induce bijections between the set of $\mu$--bitabloids and
  $\Bit$. 
\end{lemma}\Medskip

The maps defined in \eqref{maps} can be used to define an action of
$W_r$ on the set $\Bit$ by taking the preimage under one of these
maps, acting on the resulting bitabloid and taking the image under the
same map again. The action of $\sym r$ on $\Bit$ is given by
place permutations, where the numbering of the places is determined by 
$\t^\lambda$ if we use $\alpha$, by $\hat{\t}^\lambda$ if we use $\beta$,
and by $\hat{\t}_\lambda$ if we use $\gamma$. 

Our main aim in this section is to prove that the Weyl module
$W^\lambda$ has a basis which is indexed by the semistandard
$\lambda$--bitableaux of various types. Indeed, we shall show that the
weight space $W^\lambda\cap\Mmu$ (see Definition \ref{weightspacedef}) is free 
as an $R$--module with basis indexed by $\Bitss$. We begin with a special 
case.
 
\begin{thm}\label{basisgroupcase} Let $R$ be a field of characteristic
  zero and $Q=q=1$. Then 
$$
\dim(W^\lambda\cap\Mmu) = |\Bitss|.
$$
\end{thm}
\pf In this case $\H$ is isomorphic to the group algebra $RW_r$ so
$ h_{a,r-a} = w_{a,r-a}$, and
$$
\uap = \prod_{i=1}^a(1+(i,-i)) \quad \text{and}\quad 
u_{r-a}^-= \prod_{i=1}^{r-a}(1-(i,-i)).
$$

Since $R$ is a field of characteristic zero, $W^\lambda\cap\Mmu$ is
spanned by the elements of the form $\toid{\t}z_\lambda$ where $\toid{\t}$
varies over the $\mu$--bitabloids. This follows immediately from
Definition \ref{weightspacedef} and Lemma \ref{gewichtszerl}.

Let $\t$ be a row standard $\mu$--bitableau. For $i\in \rplus$, we have
$$ \toid{\t}(i,-i) = \begin{cases}
          \toid{\t}, & \text{ if $i$ is an entry of $\ta{1}$},\\
          \toid{\t(i,-i)}, & \text{ if $i$ or $-i$ is an entry of $\ta{2}$}.
\end{cases}$$
Let $X = \{i\in\rplus\strich i\text{ or }-i\in\ta{2}\}$ and 
$X_a  = X\cap \{1,2,\ldots ,a\}$. Then
\begin{equation}
\toid{\t}\uap = 2^{a-|X_a |}\sum\{\toid{\t_1}\strich \t_1\in A_\t \},
\end{equation}
where $A_\t $ is the set of row standard $\mu$--bitableaux which agree
with $\t$, except that the integers in $X_a $ are allowed to have either
sign.

Next,
\begin{equation*}
\tilde{u}_{r-a}^- :=
w_{a,r-a}u_{r-a}^-w_{a,r-a}^{-1} = \prod_{i=a+1}^r(1-(i,-i)).
\end{equation*}
Therefore
\begin{equation}
\toid{\t}\uap \tilde{u}_{r-a}^- = 
\begin{cases} 
   2^{a-|X_a |}\sum\{\pm\toid{\t_1} \strich \t_1\in B_\t \}, 
            &\text{ if }\{a+1,\ldots,r\}\subseteq X,\\
   0,  & \text{ otherwise.}
\end{cases}\end{equation}
where the sum is over the set $B_\t $ of the $2^{|X|}$ row standard
$\mu$--bitableaux $\t_1 = (\ta{1}_1,\ta{2}_1)$ which are the same as $\t$
except that the integers in $\ta{2}_1$ are allowed to have either
sign. (The coefficient of $\toid{\t_1}$ is $+1$ if and only if 
$|\{i\in\ta{2}\strich a+1\leq i\leq r\}| - 
|\{i\in\ta{2}_1\strich a+1\leq i\leq r\}|$ is even.)

Assume now that $\{a+1,\ldots,r\}\subseteq X$. Then 
$$
\toid{\t}\uap w_{a,r-a}u_{r-a}^- =
2^{a-|X_a |}\sum\{\pm\toid{\t_1}\strich \t_1\in C_\t \},
$$
where the sum is over the set $C_\t $ of row standard $\mu$--bitableaux
$\t_1 = (\ta{1}_1,\ta{2}_1)$ which are row equivalent to $\t w_{a,r-a}$,
except that the integers in $\ta{2}_1$ can have either sign. Let
$\toid{\t^+}$ be the $\mu$--bitabloid obtained from $\t w_{a,r-a}$ by
changing the signs of all the negative entries and taking the
row equivalence class. Since $\{a+1,\ldots,r\}\subseteq X$, we have
$$
\{a+1,\ldots,r\}w_{a,r-a} = \{1,2,\ldots,r-a\}\subseteq \t^{+(2)}.
$$
Hence $ \beta(\toid{\t^+}) \in \Bitplus$ and 
$ 2^{|X_a |-a}\beta(\toid{\t}\uap w_{a,r-a}u_{r-a}^-) =
\pm\sigma(\beta\toid{\t^+}). 
$
Therefore, by \eqref{2-3},
$$
2^{|X_a |-a}\beta(z_\lambda) =
\pm\sigma(\beta(\toid{\t^+})\hat{x}_\lambda
T_{\hat\pi_\lambda}\hat{y}_\lambda).
$$
If $\beta(\toid{\t^+}) = \T$, then $\T\hat{x}_\lambda$ is a multiple of the
sum over those $\T'$ which are row equivalent to $\T$ (since
$\hat{x}_\lambda$ is the sum over the row symmetrizer of $\T$). Hence
$\beta(z_\lambda)$ is a non-zero multiple of $\sigma(E_\T)$. The 
theorem now follows from Corollary \ref{domsignbase}.
\epf

We shall generalise Theorem \ref{basisgroupcase}, making use of the
special case in the course of the proof. First we need to
reformulate some results from section 4 in the language of
$\lambda$--tableaux using the map $\alpha$ introduced in
Definition~\ref{maps}. Define $\tilde{\mu}$
to be the composition of $r$ into $2n$ parts
$\tilde{\mu} = (\tilmu{1},\tilmu{2},\ldots,\tilmu{2n})$ as
follows. Suppose that $\mue 1 = (\mu^{(1)}_1,\ldots,\mu^{(1)}_{n_1})$ and 
$\mue 2 = (\mu^{(2)}_1,\ldots,\mu^{(2)}_{n_2})$. Thus $n_1+n_2 = n$. We 
define
\begin{equation}
\tilmu{i} = 
\begin{cases} 
\mu^{(1)}_i, & \text{ for } 1\leq i\leq n_1,\\
\mu^{(2)}_i, & \text{ for } n+1\leq i\leq n+n_2,\\
0,        & \text{ otherwise}.
\end{cases}\end{equation}  
Define $\tilde\lambda$ similarly.
 
Denote the set of $\tilde{\lambda}$--tableaux of type $\tilde{\mu}$
by $\Bittilde$. We can turn $\T\in\Bitplus$ 
into an element $\tilde{\T}$ of $\Bittilde$ by combining the two components 
of~$\T$.

\begin{ex}\label{beisp5} If $r=n=7$ and $\lambda,\mu\in\bicomp{n,r}$
  as in Example {\rm{\ref{beisp1}}}. Then
\begin{align*}
    \tilde{\lambda} &= (3,2,0,0,0,0,0,1,1,0,0,0,0,0)\qquad\text{and}\\
    \tilde{\mu}    &= (3,1,0,0,0,0,0,2,1,0,0,0,0,0). 
\end{align*}
Moreover if 
$
\T = \left(\begin{matrix} 1&1&1\\
                     2&8&\     
\end{matrix}\right.
      \quad , \quad
   \left.\begin{matrix}  8\\
                     9 
                                        \end{matrix}
            \right ) \in \Bitplus,
$
then
$$
\tilde{\T}\qquad = \qquad\begin{matrix} 1&1&1\\  
                       2&8&\ \\
                       -&\ &\ \\
                       -&\ &\ \\
                       -&\ &\ \\
                       -&\ &\ \\
                       -&\ &\ \\
                       8&\ &\ \\
                       9&\ &\ \\
                       -&\ &\ \\
                       -&\ &\ \\
                       -&\ &\ \\
                       -&\ &\ \\
                       -&\ &\   \end{matrix}
$$
in $\Bittilde$.
\end{ex}\Medskip

As in \eqref{situation} we let $I$ and $J$ be the subsets of $\Delta$ such
that $x_\lambda=x_J$ and $x_\mu=x_I$.

We have the following lemma. The first part is trivial; the second
part follows from part one and
the corresponding result for type $\bf A$ \cite[1.7(i)]{DJ1}. 

\begin{lemma}\label{tblident}\begin{enumerate}
\item The map $\T\mapsto\tilde{\T}$ defines a bijection between 
$\Bitplus$ and $\Bittilde$.
\item There is a bijection between $\D_{I,J}^\Delta$ and the set of
  positive row standard $\lambda$--bitableaux of type $\mu$
  given by $c\mapsto \T_c$ for $c\in\D_{I,J}^\Delta$, where 
$\T_c = \alpha\toid{\t^\mu c}$. 
\end{enumerate}\end{lemma}\Medskip

\begin{ddef}\label{rowequiv} Given $\T_1,\T_2\in\Bit$ we write $\T_1\sim
  \T_2$ if $\T_1$ and $\T_2$ are row equivalent.
\end{ddef}\Medskip

Corollary \ref{blackd} (see also \ref{tableaux sums}), taken in conjunction 
with \cite[1.7, 3.4]{DJ1} now gives the following theorem
(note that the map $A\mapsto 1_A$ in \cite[1.7]{DJ1} is the 
analogue of the map $\alpha$ from Definition~\ref{maps}).

\begin{thm}\label{genmapstotbl} Let $c\in\D_{I,J}^\Delta$. Then 
$$
\lmcphi(\toid{\tlam}) \equiv \sum_{\alpha\toid{\t_1}\sim
  \T_c}\toid{\t_1} \pmod{\Mmu_-}
$$
where $\T_c = \alpha\toid{\t^\mu c}$.
\end{thm}\Medskip

We are now prepared to embark on the generalisation of Theorem
\ref{basisgroupcase}. 

Let $\T\in\Bitss$ and let $\t$ be the row standard $\mu$--bitableau
such that $\alpha\toid{\t} = \T$, which is given by
Lemma~\ref{inversemaps}.

The semistandard tableau $\T$ corresponds to an element
$c\in\D_{I,J}^\Delta$ by Lemma~\ref{tblident}.  By 
Theorem~\ref{genmapstotbl} the $\H$-homomorphism
$\lmcphi$ maps the generator $\toid{\tlam}$ of $M^\lambda$ 
to an element $v$ of $\Mmu$ which satisfies
$$
v\equiv\sum_{\alpha\toid{\t_1}\sim \T}\toid{\t_1}\pmod{\Mmu_-}.
$$    
Since $h_{a,r-a}\in\H(\sym{r})$, 
$$
vh_{a,r-a}\equiv\sum_{\alpha\toid{\t_1}\sim \T}\toid{\t_1}h_{a,r-a}
\pmod{\Mmu_-}
$$
by \eqref{3.12}. Next, Corollary~\ref{3.15} gives
$$
vh_{a,r-a}\equiv \sum_{\alpha\toid{\t_1}\sim
  \T}r_{\t_1}\toid{\t_1w_{r,r-a}}+v_1\pmod{\Mmu_-} ,
$$
where each $r_{\t_1}$ is a unit in $R$ and $v_1$ is a linear
combination of $\mu$--bitabloids $\toid{\t_1'}$ satisfying
$ \beta\toid{\t_1'}>\beta\toid{\t w_{a,r-a}}$ 
(see part (ii) of Definition \ref{ordnung}).

Now, $\alpha\toid{\t_1}=\beta\{\t_1 w_{a,r-a}\}$, so 
\begin{equation}
vh_{a,r-a}\equiv \sum_{\beta\toid{\t_2}\sim
  \T}r_{\t_2}\toid{\t_2}+v_1\pmod{\Mmu_-},
\label{6-5}\end{equation}
where each $r_{\t_2}$ is a unit.

Note that the numbers $1,2,\ldots,r-a$ belong to
$\hat{\t}^{\lambda(2)}$. On the other hand, since $\T$ is semistandard
by assumption, all the numbers in $\T^{(2)}$ are greater than~$n$. Thus for 
all $\beta\{\t_2\}$ which are row equivalent to $\T$, all of the entries in 
the second component of $\beta\{\t_2\}$ are greater than $n$. By
Definition \ref{maps} we conclude that 
all the $\mu$--bitabloids $\toid{\t_2}$ which appear in \eqref{6-5}
have $1,2,\ldots,r-a$ in $\ta{2}$. A similar result applies to the 
$\mu$--bitabloids which occur in $v_1$. Therefore, by Corollary \ref{3.13},
\begin{equation}
vh_{a,r-a}u^-_{r-a} \equiv\sum_{\beta\toid{\t_2}\sim
  \T}r'_{\t_2}\toid{\t_2}+v_2\pmod{\Mmu_-}.
\end{equation}
where each $r'_{\t_2}$ is a unit and $v_2$ is a linear
combination of $\mu$--bitabloids $\toid{\t'_1}$ which are precisely the
$\mu$--bitabloids which are involved with nonzero coefficients in
$v_1$. (Indeed Corollary \ref{3.13} implies that the coefficients of
$\toid{\t'_1}$ in $v_1$ and $v_2$ differ only by a unit.) 

Next, let $\toid{\t^*}$ be the $\mu$--bitabloid such that
$\gamma\toid{\t^*} = \T$. Note that the row standard $\mu$--bitableau 
$\t^*$ in $\toid{\t^*}$ is $\t\hat{\pi}_\lambda$.

We have that
$$ vh_{a,r-a}u^-_{r-a}T_{\hat{\pi}_\lambda}
       \equiv r\toid{\t^*}+v_3\pmod{\Mmu_-}  $$
where $r$ is a unit and $v_3$ is a linear combination of
$\mu$--bitabloids $\toid{\t_3}$ such that $\gamma\toid{\t_3}>\T$.
To see this, note part (iii) of Definition \ref{ordnung} and compare with
\cite[7.26]{DJ5}. The matrix $\chi(\t,\t^\mu w_\mu)$ which appears in
\cite[7.26]{DJ5} is defined in such a way that its $(j,k)$th entry is equal
to the number of entries less than or equal to $j$ in the first $k$ columns
of the $\mu$--tableau of type $\lambda$ obtained by replacing each entry $i$
of $\t^\mu w_\mu$ by $\row_\t(i)$.

We arrive at the following element of $W^\lambda\cap\Mmu$.
$$ \lmcphi(z_\lambda) 
      = vh_{a,r-a}u^-_{r-a}T_{\hat{\pi}_\lambda}\hat{y}_\lambda
      \equiv r\toid{\t^*} + v_4\pmod{\Mmu_-}, $$
where $v_4$ is a linear combination of $\mu$--bitabloids $\toid{\t_4}$
such that 
$$
\gamma\toid{\t_4}>\T,
$$
by part (i) of Definition \ref{ordnung}. This is justified as follows. All of
the $\mu$--bitabloids $\toid{\t_3}$ involved in $v_3$ when acted upon by 
the terms $T_w$, appearing in $\hat{y}_\lambda$, are linear
combinations of $\mu$--bitabloids $\toid{\t_4}$ which are obtained from
$\toid{\t_3}$ by permuting the entries of the columns in $\t_3$. Since
$\gamma\toid{\t_3}>\T$, Definition~\ref{ordnung} shows that 
$\gamma\{\t_4\}>\T$. Also, $\toid{\t^*}\T_w$ is a linear combination of 
terms of the form  $\toid{\t^*w'}$ where $w'$ is an element in the column 
stabilizer of $\hat{\t}_\lambda$. From part (i) of Definition
\ref{ordnung} we conclude that $\gamma\toid{\t^*w'}>\T$ for $w\neq 1$ 
since $\T$ is semi-standard by assumption. We have shown the following.

\begin{lemma}\label{gewichtsbasis} Let $\T\in\Bitss$ and let 
$$
v_{\lambda,\mu}(\T) =
\lmcphi(\toid{\t})h_{a,r-a}u_{r-a}^-T_{\hat{\pi}_\lambda}\hat{y}_\lambda,
$$
where $\t$ is the unique row standard $\mu$--bitableau such that
$\alpha\toid{\t} = \T$. Let $\t^*$ be the row standard $\mu$--bitableau such
that $\gamma\toid{\t^*} = \T$. Then $v_{\lambda,\mu}(\T)\in W^\lambda\cap\Mmu$ and
$v_{\lambda,\mu}(\T)$ is congruent modulo $\Mmu_-$ to a linear combination of
$\mu$--bitabloids $\toid{\t_4}$ such that 
$$
\gamma\toid{\t_4} \geq \T
$$
and the coefficient of $\toid{\t^*}$ in $v_{\lambda,\mu}(\T)$ is invertible.
\end{lemma}\Medskip

Recall that by Lemma
\ref{gewichtszerl} $W^\lambda\cap\Mmu$ is a weight space of the 
$(Q,q)$--Weyl module $W^\lambda$.

\begin{cor}\label{basisgewicht} Let $\lambda,\mu\in\bicomp{n,r}$. Then
$\{v_{\lambda,\mu}(\T)\strich \T\in\Bitss\}$ is a linearly independent subset of the
weight space $W^\lambda\cap\Mmu$ of $W^\lambda$.
\end{cor}

We now prove our main result.

\begin{thm}[The Semistandard Basis Theorem]\hfil\newline\label{main}%
Let $\lambda,\mu\in\bicomp{n,r}$ and let $Q,q$ be
invertible elements of $R$. Then the weight space
$W^\lambda\cap\Mmu$ of the Weyl module $W^\lambda$
is free as an $R$--module with basis 
$\{v_{\lambda,\mu}(\T)\strich \T\in\Bitss\}$. 
Consequently, $W^\lambda$ is free as an $R$--module with basis
$$
\{v_{\lambda,\mu}(\T)\strich \mu\in\bicomp{n,r}, \T\in\Bitss\}.
$$  
\end{thm}
\pf Suppose, for the moment, that $R = {\Bbb Q}(q,Q)$ where $Q$ and
$q$ are independent transcendentals. Then $\H_{R,q,Q}(W_r)$ is 
isomorphic to $RW_r$ by \cite{BC}. Hence, in this case,
$\dim(W^\lambda\cap\Mmu) = |\Bitss|$ by Theorem \ref{basisgroupcase},
so $\{v_{\lambda,\mu}(\T)\strich \T\in\Bitss\}$ is a basis of $W^\lambda\cap\Mmu$.

Assume that $m$ is a nonzero element of $W^\lambda\cap\Mmu$ and that
the coefficient of every $\mu$--bitabloid which is involved in $m$
belongs to ${\Bbb Z}[q,q^{-1},Q,Q^{-1}]$. We write
$$
m = \sum_{\T\in\Bitss}r_\T v_{\lambda,\mu}(\T)
$$
with coefficients $r_\T$ in ${\Bbb Q}(q,Q)$. We claim that $r_\T\in{\Bbb
  Z}[q,q^{-1},Q,Q^{-1}]$. In the total order of \eqref{ordnung} let $\T_1$ be 
the first element of $\Bitss$ such that $r_{\T_1}\neq 0$,
and let $\t_1$ be the row standard $\mu$--bitableau with $\gamma\toid{\t_1}
= \T_1$. The coefficient of $\toid{\t_1}$ in $m$ belongs to 
${\Bbb Z}[q,q^{-1},Q,Q^{-1}]$, so by Lemma \ref{gewichtsbasis}
$$
r_{\T_1}\in{\Bbb Z}[q,q^{-1},Q,Q^{-1}]. 
$$
Using Lemma \ref{gewichtsbasis}, we see that
$$
m-r_{\T_1}v_{\lambda,\mu}(\T_1) = \sum_{\stackrel{\T\in\Bitss}{\T\ne\T_1}}r_\T v_{\lambda,\mu}(\T)
$$
has the property that the first element $\T_2$ of $\Bitss$ such that
$r_{\T_2}\neq 0$ satisfies $\T_2>\T_1$. Hence, by induction, every $r_\T$
belongs to ${\Bbb Z}[q,q^{-1},Q,Q^{-1}]$.

Now assume that $R$ is an arbitrary commutative ring and $q$, $Q$ are
invertible elements of $R$. The result of the last paragraph shows
that every non-zero element of $W^\lambda\cap\Mmu$ can be written as a
linear combination of $\{v_{\lambda,\mu}(\T)\strich \T\in\Bitss\}$. Taken in
conjunction with \eqref{basisgewicht}, this proves that
$$
\{v_{\lambda,\mu}(\T)\strich \T\in\Bitss\}
$$
is a basis of $W^\lambda\cap\Mmu$ and concludes the proof of the
theorem.
\epf

\it
\begin{tabular}{ll}
E--mail addresses:& rdipper@@mathematik.uni--stuttgart.de\\ 
\              & g.james@@ic.ac.uk\\
\              & a.mathas@@ic.ac.uk  
\end{tabular}
\rm
\end{document}